\documentclass[11pt,reqno]{amsart}

\usepackage[T1]{fontenc}
\usepackage{lmodern}
\usepackage{microtype}
\usepackage{amsmath,amssymb,mathtools}
\numberwithin{equation}{section}
\usepackage[all]{xy}
\usepackage{enumitem}
\usepackage{tikz-cd}
\usepackage[colorlinks=true,linkcolor=blue,citecolor=blue,urlcolor=blue]{hyperref}
\usepackage{extarrows}
\usepackage{tabularx}
\usepackage{array}

\setlist[enumerate,1]{label=(\arabic*),leftmargin=2.2em}
\setlist[itemize]{leftmargin=2em}
\newtheorem{theorem}{Theorem}[section]
\newtheorem{lemma}[theorem]{Lemma}
\newtheorem{proposition}[theorem]{Proposition}
\newtheorem{corollary}[theorem]{Corollary}

\newtheorem{remark}[theorem]{Remark}
\newtheorem{definition}[theorem]{Definition}

\newcommand{\cA}{\mathcal A}
\newcommand{\cB}{\mathcal B}
\newcommand{\cX}{\mathcal X}
\newcommand{\cY}{\mathcal Y}

\newcommand{\cC}{\mathcal C}

\newcommand{\Ext}{\operatorname{Ext}}
\newcommand{\Hom}{\operatorname{Hom}}
\newcommand{\Tor}{\operatorname{Tor}}
\newcommand{\Coker}{\operatorname{Coker}}
\newcommand{\Ker}{\operatorname{Ker}}
\newcommand{\Mod}{\text{-}\mathrm{Mod}}

\newcommand{\Com}{\operatorname{Com}}
\newcommand{\CoCom}{\operatorname{CoCom}}
\newcommand{\Ima}{\operatorname{Im}}

\newcommand{\Id}{\operatorname{Id}}
\newcommand{\radf}{\mathfrak r}
\newcommand{\coradf}{\mathfrak q}
\newcommand{\Ann}{\operatorname{Ann}}

\title[Cotorsion Pairs in $\eta$-Extensions]{Cotorsion Pairs in $\eta$-Extensions with Applications to $\theta$-Extension Rings and Tensor Rings}
\author{Pan Chen}
\address{Pan Chen: School of Mathematics,
Yunnan Normal University, Kunming, Yunnan 650500, China.}
\email{chenpan901@163.com}
\author{Yongyun Qin}
\address{Yongyun Qin: School of Mathematics,
Yunnan Normal University, Kunming, Yunnan 650500, China.}
\email{qinyongyun2006@126.com}
\author{Jinrui Yang}
\address{Jinrui Yang: School of Mathematics and statistics,
Guizhou University, Guiyang, Guizhou 100190, China.}
\email{15906981285@163.com}

\begin{document}

\begin{abstract}
Let $\cB$ be an abelian category with enough projective and injective objects, and let $\cB\ltimes_\eta F$ be an $\eta$-extension of $\cB$ with $\eta$ associative.
Given a cotorsion pair $(\cX, \cY)$ in $\cB$, we show that the cotorsion pair $(\Delta(\cX),\Delta(\cX)^\perp)$ exists for arbitrary associative $\eta$, removing the square-zero assumption $F^2=0$ from Hu's construction. When $(\cX, \cY)$ is hereditary, the heredity of $(\Delta(\cX),\Delta(\cX)^\perp)$ is characterized by the first left derived functor of a canonical defect functor $\delta$. We also remark that Hu's construction and hereditary result for the cotorsion pair
$({}^\perp U^{-1}(\cY), U^{-1}(\cY))$ remain valid without the nilpotence assumption on $\eta$. Under nilpotence hypotheses, we establish a completeness criterion for
$(\Delta(\cX),\Delta(\cX)^\perp)$ using the Beligiannis radical functor $\radf$. And we also provide a completeness criterion for $({}^\perp U^{-1}(\cY), U^{-1}(\cY))$ if
$\cB\ltimes_\eta F$ has enough injective objects. We give sufficient conditions for their coincidence. Dual results are obtained for $\zeta$-coextensions. Finally, we specialize the theory to $\theta$-extension rings and tensor rings, where the abstract conditions admit explicit module-theoretic descriptions.
\end{abstract}

\maketitle

	
	\medskip
	

\section{Introduction}\label{section:introduction}

Cotorsion pairs, introduced by Salce \cite{Salce}, are Ext-orthogonal analogues of torsion pairs. They are standard tools in relative homological algebra and approximation theory, and complete cotorsion pairs are closely related to abelian model structures through Hovey's correspondence; see, for example, \cite{EnochsJenda,Gillespie,GobelTrlifaj,Hovey}. They also provide useful tools in the study of finitistic dimension conjectures \cite{AngeleriHernandez,AngeleriTrlifaj,WangLiHu}. A recurring problem is to construct cotorsion pairs in a new category from known cotorsion pairs in other categories. This problem has been studied for triangular matrix rings, trivial extensions, Morita context rings, and recollements of abelian categories; see \cite{Mao,ZhuPengDing,MaoTrivial,CuiRongZhang,YangQin}.

The categorical framework used here is the theory of extensions of abelian categories developed by Marmaridis and Beligiannis \cite{Marmaridis,Beligiannis}. Let $\cB$ be an abelian category with enough projective and injective objects. Let $F\colon\cB\to\cB$ be a right exact functor, and let $\eta\colon F^2\to F$ be an associative natural transformation. The corresponding $\eta$-extension $\cB\ltimes_\eta F$ consists of pairs $(X,f)$ with $X$ an object of $\cB$ and $f\colon F(X)\to X$ satisfying $f\circ F(f)=f\circ \eta_X$. When $\eta=0$, the $\eta$-extension reduces to the right trivial extension. There are two natural adjoint pairs of functors connecting $\cB$ and $\cB\ltimes_\eta F$: $(T,U)$ and $(C,Z)$. See Subsection~2.2 for more on $\eta$-extensions.

Our starting point is the recent work of Hu \cite{Hu}. Hu works with $\eta$-extensions under the standing assumption that $\eta$ is associative and nilpotent. Given a cotorsion pair $(\cX,\cY)$ in $\cB$, he proves, under the condition $\mathbb{L}_1F(\cX)=0$, that
$({}^\perp U^{-1}(\cY),U^{-1}(\cY))$
is a cotorsion pair and that its heredity is equivalent to the heredity of $(\cX,\cY)$; see \cite[Theorem~3.5]{Hu}. Hu also considers a class $\Delta(\cX)$ defined by exactness of the complex
$F^2(X)\xrightarrow{\eta_X-F(f)}F(X)\xrightarrow{f}X$
together with $\Coker f\in\cX$. However, the corresponding cotorsion pair
$
(\Delta(\cX),\Delta(\cX)^\perp)
$
is constructed there only in the square-zero case $F^2=0$.  Moreover, under the additional condition $\mathbb{L}_1F(\cX)=0$, its heredity is equivalent to that of $(\cX,\cY)$; see \cite[Theorem~3.7]{Hu}. The completeness results are asymmetric as well: Hu proves completeness for the $U^{-1}$-pair in the right trivial extension case under the additional condition $F(\cX)\subseteq\cY$ \cite[Theorem~3.13]{Hu}, whereas completeness of the $\Delta$-pair is not studied separately there.

The aim of this paper is to extend the main results of \cite{Hu} to the general associative $\eta$-extension, and, at the same time, to isolate precisely where nilpotence is actually needed. Our first main result concerns the existence and heredity of the two cotorsion pairs. The exact sequence \eqref{eq:TU-seq} induced by the adjoint pair $(T,U)$
and the defect functor $\delta$ play key roles, both of which are defined in Subsection~2.3.

\medskip
\noindent\textbf{Theorem A.}
(See Theorems~\ref{thm:Delta-cotorsion}, \ref{thm:Delta-heredity} and \ref{thm:U-cotorsion}; compare \cite[Theorem~A]{Hu})
Let $\cB$ be an abelian category with enough projective and injective objects and let $\cB\ltimes_\eta F$ be an $\eta$-extension of $\cB$ with $\eta$ associative. Assume that $(\cX,\cY)$ is a cotorsion pair in $\cB$.
\begin{enumerate}
\item The pair
$
(\Delta(\cX),\Delta(\cX)^\perp)
$
is a cotorsion pair in $\cB\ltimes_\eta F$. If $(\cX,\cY)$ is hereditary, then this pair is hereditary if and only if
$
\mathbb{L}_1\delta(\Delta(\cX))=0.
$
\item If $\mathbb{L}_1F(\cX)=0$, then
$
({}^\perp U^{-1}(\cY),U^{-1}(\cY))
$
is a cotorsion pair in $\cB\ltimes_\eta F$, and it is hereditary if and only if $(\cX,\cY)$ is hereditary.
\end{enumerate}
\medskip

Part~(1) removes the square-zero restriction from Hu's $\Delta$-construction and replaces the trivial-extension hereditary obstruction $\mathbb{L}_1F$ by  $\mathbb{L}_1\delta$. When $\eta=0$, the condition $\mathbb{L}_1\delta(\Delta(\cX))=0$ is equivalent to $\mathbb{L}_1F(\cX)=0$, so Hu's hereditary criterion is recovered. Part~(2) is not new. We list it here only to show that the nilpotence assumption in Hu's $U^{-1}$-construction and hereditary theorem is unnecessary. Thus, for both cotorsion pairs, neither the existence nor the hereditary theory requires the nilpotence of $\eta$.

Nilpotence enters at the next stage. Its role is not to produce the two cotorsion pairs, but to make the Beligiannis radical functor $\radf=\Ker(\Id\to ZC)$, defined in Subsection~2.3, nilpotent and hence to provide every object $(X,f)$ with a finite $Z(\cB)$-filtration; see Lemma~\ref{lem:nilp-Z-filt}. This yields our second main result.

\medskip
\noindent\textbf{Theorem B.}
(See Theorems~\ref{thm:completeness} and \ref{thm:U-completeness}; compare \cite[Theorem~B]{Hu})
Let $\cB$ be an abelian category with enough projective and injective objects and let $\cB\ltimes_\eta F$ be an $\eta$-extension of $\cB$ with $\eta$ associative. Assume that $(\cX,\cY)$ is a complete cotorsion pair in $\cB$.
\begin{enumerate}
\item If $\eta$ is $N$-nilpotent, $(\cX,\cY)$ is hereditary,
$
\mathbb{L}_1\delta(\Delta(\cX))=0$ and
$C\radf^nT(\cX)\subseteq\cY$ for all $1\le n\le N,
$
then $(\Delta(\cX),\Delta(\cX)^\perp)$ is a hereditary complete cotorsion pair.
\item If $\cB\ltimes_\eta F$ has enough injective objects,
$
\mathbb{L}_1F(\cX)=0
$ and
$
F(\cX)\subseteq \cY,
$
then $({}^\perp U^{-1}(\cY),U^{-1}(\cY))$ is a complete cotorsion pair.
\end{enumerate}
\medskip

Completeness of the $\Delta$-pair in part~(1) is new: this problem is not addressed separately in \cite{Hu}. Part~(2), on the other hand, extends Hu's completeness theorem \cite[Theorem~3.13]{Hu} from right trivial extensions to general associative $\eta$-extensions. We remark that \cite{LG} also gives a completeness criterion for $({}^\perp U^{-1}(\cY),U^{-1}(\cY))$, which need certain nilpotence condition; see Theorem~\ref{thm:pre-U-completeness}. Our completeness criterion does not require the nilpotence condition, although we do need the category
$\cB\ltimes_\eta F$ to have enough injective objects. This additional assumption is automatically satisfied when applied to module categories.

For convenience, Table~\ref{tab:main-results} summarizes the additional conditions under which the two constructions preserve the property of being a cotorsion pair, a hereditary cotorsion pair, or a hereditary complete cotorsion pair when passing from the base category $\cB$ to the $\eta$-extension
$\cB\ltimes_\eta F$.
More precisely, the last three columns record the additional conditions required for the corresponding preservation property.
For instance, if $(\mathcal X,\mathcal Y)$ is a hereditary cotorsion pair in
$\mathcal B$, then the first construction $(\Delta(\mathcal X),\Delta(\mathcal X)^\perp)$ yields a hereditary
cotorsion pair in $\cB\ltimes_\eta F$ under the condition $\mathbb L_1\delta(\Delta(\cX))=0$, while the second
construction
$({}^\perp U^{-1}(\cY),U^{-1}(\cY))$ yields one provided
that $\mathbb L_1F(\cX)=0$.

\begin{table}[htbp]
\centering
\caption{Preservation of cotorsion-pair properties under the two constructions}
\label{tab:main-results}
\renewcommand{\arraystretch}{1.3}
\setlength{\tabcolsep}{5pt}

\begin{tabularx}{\textwidth}{
|>{\raggedright\arraybackslash}p{0.25\textwidth}
|>{\raggedright\arraybackslash}p{0.18\textwidth}
|>{\raggedright\arraybackslash}p{0.22\textwidth}
|>{\raggedright\arraybackslash}X|}
\hline
Construction
& Cotorsion pair
& Hereditary cotorsion pair
& Hereditary complete cotorsion pair \\
\hline

$\bigl(\Delta(\mathcal X),\Delta(\mathcal X)^\perp\bigr)$
&
No additional condition.
&
$\mathbb L_1\delta(\Delta(\mathcal X))=0$.
&
\begin{minipage}[t]{\linewidth}
\raggedright
$\mathbb L_1\delta(\Delta(\mathcal X))=0$;
\medskip
$\eta$ is $N$-nilpotent;
\medskip
$C\radf^nT(\mathcal X)\subseteq\mathcal Y$
for $1\leq n\leq N$.
\vspace{2pt}
\end{minipage}
\\
\hline

$\bigl({}^\perp U^{-1}(\mathcal Y),U^{-1}(\mathcal Y)\bigr)$
&
$\mathbb L_1F(\mathcal X)=0$.
&
$\mathbb L_1F(\mathcal X)=0$.
&
\begin{minipage}[t]{\linewidth}
\raggedright
$\mathbb L_1F(\mathcal X)=0$;

\medskip
$F(\mathcal X)\subseteq\mathcal Y$;

\medskip
either $\mathcal B\ltimes_\eta F$
has enough injective objects or $\eta$ is $N$-nilpotent.
\end{minipage}
\\
\hline
\end{tabularx}
\end{table}

In addition to studying the two cotorsion pairs in parallel, we shall compare them. Under the condition $\mathbb{L}_1F(\cX)=0$, we always have
$\Delta(\cX)^\perp\subseteq U^{-1}(\cY),$
or equivalently,
${}^\perp U^{-1}(\cY)\subseteq\Delta(\cX).$
Under the nilpotence hypothesis we give sufficient conditions for these inclusions to become equalities. The key point is a splitting criterion for induced objects: for $A\in\Delta(\cX)$, $A\cong T(C(A))$ if and only if the following canonical sequence
\[
0\longrightarrow\delta(A)\longrightarrow U(A)\longrightarrow C(A)\longrightarrow0
\]
is split exact; see Lemma~\ref{lem:T-equi}.
This leads to a sufficient condition, $\Ext^1_{\cB}(C(A),\delta(A))=0\,(A\in \Delta(\cX))$, for $\Delta(\cX)=T(\cX)$ and hence for the coincidence of the two cotorsion pairs.
When $\eta=0$, this sufficient condition reduces to the given one, $\Ext^1_{\cB}(X,F(X))=0\, (X\in \cX)$, in \cite[Lemma~3.12]{Hu}. See Theorem~\ref{thm:two-equal} and Remark~\ref{rem:equal-reduce}.

The theory has a genuine dual. Given a left exact endofunctor $G$ and a coassociative natural transformation $\zeta\colon G\to G^2$, we have a
$\zeta$-coextension $G\rtimes_{\zeta}\cB$. For it, we develop the corresponding class $\nabla(\cY)$, the codefect functor $\delta^\vee$, and the coradical functor $\coradf$.
Given a cotorsion pair $(\cX,\cY)$ in $\cB$, we have two cotorsion pairs $({}^\perp\nabla(\cY), \nabla(\cY))$ and $(U^{-1}(\cX), U^{-1}(\cX)^\perp)$ in $G\rtimes_{\zeta}\cB$.
The hereditary obstruction of $({}^\perp\nabla(\cY), \nabla(\cY))$ is $\mathbb{R}^1\delta^\vee$, and the nilpotence of $\zeta$ yields the dual completeness theorem.
For the category $\cB\ltimes_\eta F$, when $F$ admits a right adjoint $G$, it may be viewed as a $\zeta$-coextension; consequently a single cotorsion pair in $\cB$ gives four natural cotorsion pairs in the $\eta$-extension $\cB\ltimes_\eta F$. See Subsections~3.5 and 3.6.

We then specialize the categorical theory to modules over rings. For a $\theta$-extension ring $\Lambda=R\ltimes_\theta L$, where $R$ is a unital ring, $L$ is an $R$-bimodule and
$\theta\colon L\otimes_R L\to L$ is an associative $R$-bimodule homomorphism,
the abstract defect functor and codefect functor become
\[
\delta(A)\cong L\otimes_\Lambda A,
\qquad
\delta^\vee(A)\cong\Hom_\Lambda(L,A),
\]
while the Beligiannis radical functor and coradical functor are described by
\[
\radf(A)\cong LA,
\qquad
\coradf(A)\cong A/\Ann_A(L).
\]
See Propositions~\ref{prop:theta-functors} and \ref{prop:theta-layers}. Thus the abstract heredity and completeness criteria become ordinary Tor-Ext conditions together with conditions on certain quotient modules.

Finally, we consider the case of the tensor rings. Let $M$ be an $R$-bimodule. Denote $M^{0}=R$ and $M^{i+1}=M\otimes_R M^{i}$ for $i\geq 0$.
The tensor ring $T_R(M) = \bigoplus_{i\geq 0}M^i$ can be viewed as a $\theta$-extension ring with a natural associative $\theta$:
$
 T_R(M) = R\ltimes_\theta M_+,
 $
where
$
M_+=\bigoplus_{i\ge1}M^i.
$
Since the tensor ring has a grading structure, all the needed functors admit simpler descriptions in this case. In particular, we have the isomorphisms
\[
C\radf^nT(X)\cong M^n\otimes_RX,
\qquad
K\coradf^nH(Y)\cong\Hom_R(M^n,Y).
\]
See Proposition~\ref{prop:tensor-layers}. It also yields concrete sufficient conditions for the coincidence of the two induced cotorsion pairs; see Corollary~\ref{cor:tensor-coincidence}. In this way, the sections on $\theta$-extension rings and tensor rings provide, respectively, a module-theoretic realization and a graded simplification of the cotorsion theory of $\eta$-extensions.

The paper is organized as follows. Section~2 recalls some definitions concerning cotorsion pairs and
$\eta$-extensions, as well as some basic facts needed later. In Section~3, we construct and study in parallel the two cotorsion pairs in
$\eta$-extensions, examine their heredity and completeness, compare the two pairs, and develop the dual theory for $\zeta$-coextension.
Section~4 specializes the abstract functors to
$\theta$-extension rings, yielding the cotorsion theory for such rings. Section~5 deals with tensor rings, where the grading structure leads to more explicit results.

\section{Preliminaries}\label{section:pre}

In this section we fix our notations, recall some relevant concepts and give some basic facts. Throughout the paper, the base category $\cB$ is assumed to be an abelian category with enough projective and injective objects, and all classes of objects in abelian categories are assumed to be closed under isomorphisms. All rings are nonzero associative rings with identity and all modules are unitary. For a ring $R$, we adopt the convention that an $R$-module is a left $R$-module. We write $R$-Mod for the category of left $R$-modules.
${}_RM$ denotes a left $R$-module and $M_R$ denotes a right $R$-module.

\subsection{Cotorsion pairs} For a class $\cC$ of objects in an abelian category $\cA$, set
\[
{}^\perp\cC
 =\{A\in\cA\mid \Ext^1_{\cA}(A,C)=0\text{ for all }C\in\cC\},
\]
\[
\cC^\perp
 =\{A\in\cA\mid \Ext^1_{\cA}(C,A)=0\text{ for all }C\in\cC\}.
\]
We write $\Ext^1_{\cA}(A,\cC)=0$ to mean $\Ext^1_{\cA}(A,C)=0\text{ for all }C\in\cC$, and similar conventions will be used.

\begin{definition}[\cite{Salce}]
A pair $(\cX,\cY)$ of classes in an abelian category $\cA$ is a \emph{cotorsion pair} if
$
\cX={}^\perp\cY
$
and
$
\cY=\cX^\perp.
$
\end{definition}

Let $M\in \cA$ and let $\cC$ be a class of objects in $\cA$. An epimorphism $f: C \rightarrow M$ with $C\in \cC$ is called \emph{a special $\cC$-precover}
 of $M$ if $\Ker(f)\in \cC^\perp$.
If every object of $\cA$
has a special $\cC$-precover, then $\cC$ is called \emph{special precovering}. Dually, A monomorphism $g: M\rightarrow C$ with $C\in \cC$
 is called \emph{a special $\cC$-preenvelope} of  $M$ if $\Coker(g)\in {}^\perp\cC$. If every object of $\cA$
has a special $\cC$-preenvelope, then $\cC$ is called \emph{special preenveloping}. See \cite[Definition 2.1.12]{GobelTrlifaj}.

A cotorsion pair $(\cX,\cY)$ is called \emph{hereditary} if $\cX$ is closed under kernels of epimorphisms and $\cY$ is closed under cokernels of monomorphisms. We shall use the following standard one-sided characterizations.

\begin{proposition}\label{prop:one-sided-criteria}
{\rm (\cite[Theorem~1.2.10]{GarciaRozas})}
Let $(\cX,\cY)$ be a cotorsion pair in an abelian category $\cA$.
\begin{enumerate}[label=\textup{(\arabic*)}]
\item Assume that $\cA$ has enough projective objects. Then $(\cX,\cY)$ is hereditary if and only if $\cX$ is closed under kernels of epimorphisms.
\item Assume that $\cA$ has enough injective objects. Then $(\cX,\cY)$ is hereditary if and only if $\cY$ is closed under cokernels of monomorphisms.
\item Assume that $\cA$ has enough projective objects or enough injective objects. Then $(\cX,\cY)$ is hereditary if and only if $\Ext^2_\cA(\cX,\cY)=0$.
\end{enumerate}
\end{proposition}

A cotorsion pair $(\cX,\cY)$ is called \emph{complete} if $\cX$ is special precovering and $\cY$ is special preenveloping. The corresponding one-sided criteria are as follows.

\begin{proposition}\label{prop:one-sided-complete}
{\rm (\cite[Proposition~7.1.7]{EnochsJenda})}
Let $(\cX,\cY)$ be a cotorsion pair in an abelian category $\cA$.
\begin{enumerate}[label=\textup{(\arabic*)}]
\item Assume that $\cA$ has enough projective objects. Then $(\cX,\cY)$ is complete if and only if $\cY$ is special preenveloping.
\item Assume that $\cA$ has enough injective objects. Then $(\cX,\cY)$ is complete if and only if $\cX$ is special precovering.
\end{enumerate}
\end{proposition}

The following proposition will be useful in our discussion of the completeness of cotorsion pairs.

\begin{proposition}\label{prop:extension-approx}
{\rm (\cite[Theorem~3.1]{AkinciAlizade})}
Let $\cA$ be an abelian category with enough projective objects or enough injective objects, and let
$
0\to A\to B\to C\to 0
$
be exact. Suppose that $(\cX,\cY)$ is a hereditary cotorsion pair in $\cA$.
\begin{enumerate}[label=\textup{(\arabic*)}]
\item If $A$ and $C$ have special $\cX$-precovers, then $B$ has a special $\cX$-precover.
\item If $A$ and $C$ have special $\cY$-preenvelopes, then $B$ has a special $\cY$-preenvelope.
\end{enumerate}
\end{proposition}

\subsection{\texorpdfstring{$\eta$}{eta}-extensions and standard functors} Let $\cB$ be an abelian category with enough projective and injective objects, $F\colon\cB\to\cB$ a right exact functor, and $\eta\colon F^2\to F$ a natural transformation. The category $\cB\ltimes_\eta F$, called an $\eta$-extension of $\cB$, is defined as follows. The objects of $\cB\ltimes_\eta F$ are pairs $(X,f)$, where $X$ is an object of $\cB$ and $f\colon F(X)\to X$ is a morphism of $\cB$ satisfying
$f\circ F(f)=f\circ\eta_X.$
A morphism from $(X,f)$ to $(Y,g)$ is a morphism $\alpha\colon X\to Y$ of $\cB$ such that
$g\circ F(\alpha)=\alpha\circ f.$

From \cite{Marmaridis}, the $\eta$-extension $\cB\ltimes_\eta F$ is an abelian category. The sequence
$(X,f)\xrightarrow{\alpha} (Y,g)\xrightarrow{\beta} (Z,h)$ in $\cB\ltimes_\eta F$ is exact if and only if the underlying sequence
$X\xrightarrow{\alpha} Y\xrightarrow{\beta} Z$ in $\cB$ is exact. When $\eta=0$, we write $\cB\ltimes F$ and call it the right trivial extension of $\cB$ by $F$ \cite{FossumGriffithReiten}.

The natural transformation $\eta\colon F^2\to F$ is \emph{associative} if
$\eta_X\circ F(\eta_X)=\eta_X\circ\eta_{F(X)}$
for every object $X\in\cB$. It is \emph{nilpotent} if there exists $N\ge1$ such that
$\eta_X\circ\eta_{F(X)}\circ\cdots\circ\eta_{F^{N-1}(X)}=0,$
for every object $X\in\cB$. In this case, we say that $\eta$ is \emph{$N$-nilpotent}. Throughout the paper, $\eta$ is assumed to be associative.

From \cite{Beligiannis}, there are the following standard functors:
\[
\xymatrix@C=3pc{
\cB \ar@<-1ex>[r]_{Z\ \ \ } & \cB\ltimes_\eta F \ar@<-1ex>[l]_{C\ \ \ } \ar@<-1ex>[r]_{\ \ \ U}
& \cB \ar@<-1ex>[l]_{\ \ \ T}
}
\]
\begin{itemize}
 	\item $U(X,f)=X$ and $U(\alpha)=\alpha$;
 	\item $Z(X)=(X,0)$ and $Z(\alpha)=\alpha$;
 	\item $T(X)=\left(X\oplus F(X),t_X =\left(\begin{smallmatrix}
 		0&0\\
 		1&\eta_X
 	\end{smallmatrix}\right)\right)$ and $T(\alpha)=\begin{pmatrix}
 	\alpha&0\\
 	0& F(\alpha)
 	\end{pmatrix}$;
 	\item $C(X,f)=\Coker(f)$ and $C(\alpha)$ is the morphism induced by the cokernel.
\end{itemize}
Note that the associativity of $\eta$ guarantees that
$T(X)$ belongs to $\cB\ltimes_\eta F$.

We recall from \cite{Beligiannis} some facts about these functors.
The pairs $(T,U)$ and $(C,Z)$ are adjoint pairs of functors. We have the identifications: $CT=\mathrm{Id}_{\cB}$, and $CZ=\mathrm{Id}_{\cB}$.
The functors $U$ and $Z$ are exact.
If $\cB$ has enough projective objects, then so does $\cB\ltimes_\eta F$. By convention,
$\cB$ has enough projective and injective objects, and hence $\cB\ltimes_\eta F$
has enough projective objects.

\vspace{-1.3\baselineskip}
\subsection{The defect functor and the radical functor}
We examine the adjoint pairs $(T,U)$ and $(C,Z)$ in detail.
First we focus on the pair $(T,U)$. The counit $\varepsilon\colon TU\to \Id_{\cB\ltimes_\eta F}$ of the adjunction $(T,U)$ is defined by
$\varepsilon_{(X,f)}=(1,f)\colon X\oplus F(X)\to X$. Obviously $\varepsilon_{(X,f)}$ is an epimorphism and it yields an exact sequence
\begin{equation}\label{eq:pre-TU-seq}
0\rightarrow \Ker\varepsilon_{(X,f)}\xrightarrow{\omega_{(X,f)}} TU(X,f)\xrightarrow{\varepsilon_{(X,f)}} (X,f)\rightarrow 0.
\end{equation}
Following \cite[Subsection~5.1]{Beligiannis}, we define an endofunctor $\Omega\colon \cB\ltimes_\eta F\to \cB\ltimes_\eta F$ by setting, for each object
$(X,f)$, $\Omega(X,f)$ to be $\Ker\varepsilon_{(X,f)}$, and for each morphism $\alpha$,
$\Omega(\alpha)$ to be the morphism induced by the universal property of the kernel. It's easy to see that
$\Omega(X,f)=(F(X),\eta_X-F(f))$ and $\omega_{(X,f)}=\binom{f}{-1}$. Thus \eqref{eq:pre-TU-seq} becomes
\begin{equation}\label{eq:TU-seq}
0\rightarrow \Omega(X,f)\xrightarrow{\omega_{(X,f)}} TU(X,f)\xrightarrow{\varepsilon_{(X,f)}} (X,f)\rightarrow 0.
\end{equation}

Applying the right exact functor $C$ to the above exact sequence \eqref{eq:TU-seq} and noting that $CT=\Id_{\cB}$, we get an exact sequence
\begin{equation}\label{eq:pre-delta-seq}
C\Omega(X,f)\xrightarrow{C(\omega_{(X,f)})} U(X,f)\xrightarrow{C(\varepsilon_{(X,f)})} C(X,f)\rightarrow 0.
\end{equation}
Set $\gamma_{(X,f)}=C(\omega_{(X,f)})$ and $\rho_{(X,f)}=C(\varepsilon_{(X,f)})$. Obviously $\gamma$ and $\rho$ are natural transformations.
We obtain the following commutative diagram
\[
\xymatrix{
& F^2(X) \ar[r]^{\binom{F(f)}{-1}} \ar[d]^{\eta_{X}-F(f)} & F(X)\oplus F^2(X) \ar[r]^{\ \ \ \ \ (1,F(f))}\ar[d]^{t_X} & F(X)\ar[d]^{f} \ar[r] & 0\\
0\ar[r] & F(X) \ar[r]^{\binom{f}{-1}} \ar[d]^{\pi_{(X,f)}} & X\oplus F(X) \ar[r]^{\ \ \ \ (1,f)} \ar[d]^{(1,0)} & X \ar[r] \ar[d]^{\rho_{(X,f)}} & 0\\
& \Coker(\eta_X-F(f)) \ar[r]^{\ \ \ \ \ \ \ \ \gamma_{(X,f)}} & X \ar[r]^{\rho_{(X,f)}} & \Coker(f) \ar[r] & 0
}
\]
in which all rows and columns are exact. We remark that $\rho_{(X,f)}$ is just the morphism corresponding to the cokernel $\Coker(f)$ and that
$\Coker(\gamma_{(X,f)})=\Coker(f)$.

\begin{definition}
Keep the above notation. Let $\delta=C\Omega$. We call $\delta$ the defect functor. $\delta(X,f)=\Coker(\eta_X-F(f))$ and $\delta$ is right exact.
\end{definition}

With the above notation, we can rewrite the exact sequence \eqref{eq:pre-delta-seq} as
\begin{equation}\label{eq:delta-seq}
\delta(X,f)\xrightarrow{\gamma_{(X,f)}} U(X,f)\xrightarrow{\rho_{(X,f)}} C(X,f)\rightarrow 0.
\end{equation}
The first left derived functor $\mathbb{L}_1\delta$, as we shall see, measures precisely the failure of the cotorsion pair
$(\Delta(\cX),\Delta(\cX)^\perp)$ to be hereditary.

We next consider the adjoint pair $(C,Z)$. The unit $\lambda\colon \Id_{\cB\ltimes_{\eta}F}\to ZC$ is defined by
$\lambda_{(X,f)}=\rho_{(X,f)}\colon X\rightarrow \Coker(f)$. It is an epimorphism and yields a short exact sequence
\begin{equation}\label{eq:ZC-seq}
0 \rightarrow \Ker(\lambda_{(X,f)}) \xrightarrow{\iota_{(X,f)}} (X,f) \xrightarrow{\lambda_{(X,f)}} ZC(X,f)\rightarrow 0.
\end{equation}

Let $F(X)\xrightarrow{k}\Ima(f)\xrightarrow{i} X$ be the canonical factorization of $f$. Thus we have the following commutative diagram with exact rows
\[
\xymatrix{
& F(\Ima(f))\ar[r]^{F(i)}\ar[d]^{k\circ F(i)} & F(X) \ar[r]^{F(\rho_{(X,f)})\ \ \ \ }\ar[d]^{f} & F(\Coker(f)) \ar[d]^{0} \\
0 \ar[r] & \Ima(f) \ar[r]^i & X \ar[r]^{\rho_{(X,f)}\ \ \ \ } & \Coker(f) \ar[r] & 0
}
\]
It is straightforward to verify that $(\Ima(f),k\circ F(i))\in \cB\ltimes_{\eta}F$; see \cite[Lemma 3.14]{Hu}. Hence $\Ker(\lambda_{(X,f)})=(\Ima(f),k\circ F(i))$.

\begin{definition}
Keep the above notation. Define a functor $\radf\colon \cB\ltimes_{\eta}F\to \cB\ltimes_{\eta}F$ as follows. For each object $(X,f)$, $\radf(X,f)=\Ker(\lambda_{(X,f)})=(\Ima(f),k\circ F(i))$. For each morphism $\alpha$,
$\radf(\alpha)$ is the morphism induced by the kernel. We call it the {radical functor} of Beligiannis.
\end{definition}

Under this definition, \eqref{eq:ZC-seq} becomes
\begin{equation}\label{eq:radical-seq}
0 \rightarrow \radf(X,f) \xrightarrow{\iota_{(X,f)}} (X,f) \xrightarrow{\lambda_{(X,f)}} ZC(X,f)\rightarrow 0.
\end{equation}
Beligiannis introduced this endofunctor $\radf$ and showed that $\radf$ is a radical of $\cB\ltimes_\eta F$ in the sense that
$\radf(A/\radf(A))=0$ for all $A\in \cB\ltimes_\eta F$ in \cite[Subsection~5.1]{Beligiannis}. Liu and Gao also defined this functor in \cite[Proposition~3.3]{LG}.
It will play a crucial role in our investigation of the completeness of cotorsion pairs.

\subsection{Some auxiliary lemmas}

In this subsection we collect some lemmas that will be used later.

%

Inspired by the proof of \cite[Theorem~3.13]{Hu}, we have the following lemma. The short exact sequence in this lemma was also  established in \cite[Lemma~4.5]{LG}. For the reader's convenience, we provide a proof.

\begin{lemma}\label{lem:E_Y-exact}
Let $0\rightarrow B\xrightarrow{j} Y \xrightarrow{p} X\rightarrow 0$ be an exact sequence in $\cB$. Set
\[
E_Y = \left(Y\oplus F(X), \beta=\begin{pmatrix}0 & 0\\ F(p) & \eta_X\end{pmatrix}\right).
\]
Then
\[
0\rightarrow Z(B)\xrightarrow{\big(\begin{smallmatrix}j\\ 0\end{smallmatrix}\big)} E_Y\xrightarrow{\left(\begin{smallmatrix}p & 0\\ 0 & 1\end{smallmatrix}\right)} T(X)\rightarrow 0
\]
is exact in $\cB\ltimes_\eta F$. Applying the functor $C$ to this exact sequence recovers the original one in $\cB$.
\end{lemma}
\begin{proof}
We have the following commutative diagram with exact rows
\[
\xymatrix@C=2.5pc{ & F(B) \ar[r]^{\big(\begin{smallmatrix}F(j)\\ 0\end{smallmatrix}\big)\ \ \ \ } \ar[d]^{0} & F(Y)\oplus F^2(X)
\ar[r]^{\big(\begin{smallmatrix}F(p)&0\\ 0&1\end{smallmatrix}\big)} \ar[d]^{\beta}
& F(X)\oplus F^2(X) \ar[d]^{t_{X}}  \\
0\ar[r] & B\ar[r]^{\big(\begin{smallmatrix}j\\ 0\end{smallmatrix}\big)\ \ \ \ } & Y\oplus F(X)
\ar[r]^{\big(\begin{smallmatrix}p&0\\ 0&1\end{smallmatrix}\big)}& X\oplus F(X)\ar[r]&0}
\]
It is easy to check that $E_Y=(Y\oplus F(X), \beta)\in \cB\ltimes_\eta F$ and thus the sequence
\begin{equation}\label{eq:E_y}
0\rightarrow Z(B)\xrightarrow{\big(\begin{smallmatrix}j\\ 0\end{smallmatrix}\big)} E_Y\xrightarrow{\left(\begin{smallmatrix}p & 0\\ 0 & 1\end{smallmatrix}\right)} T(X)\rightarrow 0
\end{equation}
is exact in $\cB\ltimes_\eta F$. Noting that $F(p)$ is epic, by a direct computation, we have $\Coker(F(B)\xrightarrow{0} B)=B\xrightarrow{1}B, \Coker(\beta)=Y\oplus F(X)\xrightarrow{(1,0)}Y$ and $\Coker(t_X)=X\oplus F(X)\xrightarrow{(1,0)} X$. Moreover, the following diagram
\[
\xymatrix{
0\ar[r] & B\ar[r]^{\big(\begin{smallmatrix}j\\ 0\end{smallmatrix}\big)\ \ \ \ }\ar[d]^{1} & Y\oplus F(X)
\ar[r]^{\big(\begin{smallmatrix}p&0\\ 0&1\end{smallmatrix}\big)}\ar[d]^{(1,0)} & X\oplus F(X)\ar[r] \ar[d]^{(1,0)} &0\\
0\ar[r] & B\ar[r]^{j} & Y\ar[r]^{p} & X\ar[r] & 0
}
\]
is commutative. Therefore, applying the functor $C$ to the exact sequence \eqref{eq:E_y} just recovers the original one in $\cB$.
\end{proof}

The following simple fact will be used frequently in the sequel.

\begin{lemma}\label{lem:L_1-seq}
Let $\beta\colon \cA_1\to \cA_2$ be a right exact functor between abelian categories, where $\cA_1$ has enough projective objects. For $X\in \cA_1$, the following
statements are equivalent:
\begin{enumerate}[label=\textup{(\arabic*)}]
\item $\mathbb{L}_1\beta(X)=0$;
\item For any exact sequence $0\rightarrow X_1\rightarrow X_2\rightarrow X\rightarrow 0$, the induced sequence
$0\rightarrow \beta(X_1)\rightarrow \beta(X_2)\rightarrow \beta(X)\rightarrow 0$ is exact. In other words, the functor $\beta$ preserves the exactness of short exact sequences with $X$
as the right end.
\end{enumerate}
\end{lemma}
\begin{proof}
{\rm (1) $\Rightarrow$ (2):} It is straightforward.

{\rm (2) $\Rightarrow$ (1):} Since $\cA_1$ has enough projective objects, we can choose an exact sequence $0\rightarrow K\rightarrow P\rightarrow X\rightarrow 0$ with $P$ projective.
Thus we have a long exact sequence
\[
\mathbb{L}_1\beta(P)=0\rightarrow \mathbb{L}_1\beta(X)\rightarrow \beta(K)\rightarrow \beta(P)\rightarrow \beta(X)\rightarrow 0.
\]
On the other hand, by {\rm (2)} we see that $\beta(K)\rightarrow \beta(P)$ is a monomorphism. Hence $\mathbb{L}_1\beta(X)=0$.
\end{proof}

\section{Cotorsion pairs in extensions of abelian categories}\label{section:extension}

Let $\cB\ltimes_\eta F$ be an $\eta$-extension of an abelian category $\cB$. In this section, we extend the main
results of \cite{Hu} on cotorsion pairs to the associative $\eta$-extension $\cB\ltimes_\eta F$.

\subsection{Existence of the \texorpdfstring{$\Delta$}{Delta}-pair} Following \cite{Hu}, for $(X,f)\in\cB\ltimes_\eta F$, we write
\[
\Com(X,f):
F^2(X)\xrightarrow{\,\eta_X-F(f)\,}F(X)\xrightarrow{\,f\,}X.
\]
For a class $\cC\subseteq\cB$, define
\[
\Delta(\cC)
 =\{(X,f)\in \cB\ltimes_\eta F \mid \Com(X,f)\text{ is exact and }\Coker(f)\in\cC\},
\]
It is easy to check that $T(\cC)\subseteq \Delta(\cC)$.

We characterize the property that $\Com(X,f)$ is exact.

\begin{lemma}\label{lem:Com-ex}
Let $(X,f)\in \cB\ltimes_\eta F$. The following statements are equivalent:
\begin{enumerate}[label=\textup{(\arabic*)}]
\item $\Com(X,f)$ is exact;
\item $\gamma_{(X,f)}\colon \delta(X,f)=\Coker(\eta_X-F(f))\to X$ is a monomorphism;
\item $\mathbb{L}_1C(X,f)=0$.
\end{enumerate}
\end{lemma}
\begin{proof}
{\rm (1) $\Leftrightarrow$ (2):} See \cite[Lemma 3.2]{Hu}.

{\rm (2) $\Rightarrow$ (3):} For any exact sequence $0\rightarrow (X_1,f_1)\rightarrow (X_2,f_2)\rightarrow (X,f)\rightarrow 0$, noting that
$\delta$ is right exact and $\gamma$ is a natural transformation,
we have the following commutative diagram with
exact rows
\[
\xymatrix{
&\delta(X_1,f_1)\ar[r]\ar[d]^{\gamma_{(X_1,f_1)}} & \delta(X_2,f_2)\ar[r]\ar[d]^{\gamma_{(X_2,f_2)}} & \delta(X,f) \ar[r]\ar[d]^{\gamma_{(X,f)}} & 0\\
0 \ar[r] & X_1\ar[r] & X_2 \ar[r] & X \ar[r] & 0
}
\]
By the snake lemma, we obtain an exact sequence $0=\Ker(\gamma_{(X,f)})\rightarrow C(X_1,f_1)\rightarrow C(X_2,f_2)\rightarrow C(X,f)\rightarrow 0$. It follows from
Lemma~\ref{lem:L_1-seq} that $\mathbb{L}_1C(X,f)=0$.

{\rm (3) $\Rightarrow$ (2):} From \eqref{eq:TU-seq} in Subsection~2.3, we have a short exact sequence
\[
0\rightarrow \Omega(X,f) \xrightarrow{\omega_{(X,f)}} TU(X,f)\xrightarrow{\varepsilon_{(X,f)}} (X,f)\rightarrow 0.
\]
Since $\mathbb{L}_1C(X,f)=0$, applying the functor $C$ to this exact sequence gives the following exact sequence
\[
0\rightarrow \delta(X,f) \xrightarrow{\gamma_{(X,f)}} X\xrightarrow{\rho_{(X,f)}} C(X,f)\rightarrow 0.
\]
Thus $\gamma_{(X,f)}$ is a monomorphism.
\end{proof}

\begin{lemma}\label{lem:CZ-ext}
{\rm (\cite[Lemma~3.3]{Hu})}
Let $A=(X,f)\in \cB\ltimes_\eta F$, $Y\in \cB$. If $\Com(X,f)$ is exact, then $Ext^1_{\cB}(C(A),Y)\cong Ext_{\cB\ltimes_\eta F}^1(A,Z(Y))$.
\end{lemma}

When $\eta=0$ and $F^2=0$, \cite[Theorem 3.7]{Hu} shows that if $(\cX,\cY)$ is a cotorsion pair in $\cB$, then
$(\Delta(\cX),\Delta(\cX)^{\perp})$ is a cotorsion pair in $\cB\ltimes_\eta F$.
In fact, it remains valid for arbitrary associative $\eta$ and right exact $F$.

We first give a useful characterization of $\Delta(^\perp\cC)$. This generalises the result in \cite[Lemma 3.8]{Hu}, which is stated only for the special case of
a right trivial extension
$\cB\ltimes F$ with $F^2=0$.

\begin{lemma}\label{lem:Delta-orthogonal}
Let $\cC$ be a class of objects in $\cB$ containing all injective objects. Then
$
\Delta({}^\perp\cC)={}^\perp Z(\cC).
$
\end{lemma}

\begin{proof}
Assume that $(X,f)\in\Delta({}^\perp\cC)$. Then $\Com(X,f)$ is exact and $C(X,f)\in{}^\perp\cC$. Lemma~\ref{lem:CZ-ext} gives
\[
\Ext^1_{\cB\ltimes_\eta F}((X,f),Z(Y))
 \cong \Ext^1_{\cB}(C(X,f),Y)=0
\]
for every $Y\in\cC$. So $(X,f)\in{}^\perp Z(\cC)$. Hence $\Delta({}^\perp\cC)\subseteq {}^\perp Z(\cC).$

Conversely, assume that $A=(X,f)\in{}^\perp Z(\cC)$. Since $\cB\ltimes_\eta F$ has enough projective objects, we can choose an exact sequence
\begin{equation}\label{eq:proj-seq}
0\longrightarrow L\xlongrightarrow{j} P\longrightarrow A\longrightarrow 0
\end{equation}
with $P$ projective.
Since $\cB$ has enough injective objects, there exists a monomorphism $i\colon C(L)\rightarrow I$ with $I$ injective. Applying $\Hom_{\cB\ltimes_\eta F}(-, Z(I))$ to the preceding
exact sequence \eqref{eq:proj-seq}, we have the following exact sequence
\[
\Hom_{\cB\ltimes_\eta F}(P, Z(I))\xrightarrow{j^*} \Hom_{\cB\ltimes_\eta F}(L, Z(I))\to \Ext^1_{\cB\ltimes_\eta F}(A, Z(I)).
\]
Noting that $\cC$ contains all injective objects, we know that $\Ext^1_{\cB\ltimes_\eta F}(A, Z(I))=0$ from $A\in{}^\perp Z(\cC)$. Thus $j^*$ is an epimorphism.
Since $(C,Z)$ is an adjoint pair, we have a commutative diagram
\[
\xymatrix{
\Hom_{\cB}(C(P), I)\ar[r]^{C(j)^*}\ar[d]^{\cong} & \Hom_{\cB}(C(L), I)\ar[d]^{\cong}\\
\Hom_{\cB\ltimes_\eta F}(P, Z(I))\ar[r]^{j^*} & \Hom_{\cB\ltimes_\eta F}(L, Z(I))
}
\]
and hence $C(j)^*$ is an epimorphism. There is a morphism $u\in \Hom_{\cB}(C(P), I)$ such that $C(j)^*(u)=i$, i.e., $u\circ C(j)=i$. Since $i$ is a monomorphism, $C(j)$ is monic.
Applying the right exact functor $C$ to the preceding short exact sequence \eqref{eq:proj-seq}, we get a long exact sequence
\[
0=\mathbb{L}_1C(P)\to\mathbb{L}_1C(A)\to C(L) \xrightarrow{C(j)} C(P)\to C(A)\to 0.
\]
Therefore, $\mathbb{L}_1C(A)=0$. By Lemma~\ref{lem:Com-ex}, $\Com(X,f)$ is exact.

Lemma~\ref{lem:CZ-ext} yields
\[
\Ext^1_{\cB}(C(X,f),Y)
 \cong \Ext^1_{\cB\ltimes_\eta F}((X,f),Z(Y))=0
\]
for any $Y\in\cC$. We get that $\Coker(f)\in{}^\perp\cC$, and hence $(X,f)\in\Delta({}^\perp\cC)$. Therefore ${}^\perp Z(\cC)\subseteq \Delta({}^\perp\cC).$
\end{proof}

From Lemma~\ref{lem:Delta-orthogonal}, we immediately obtain the desired cotorsion pair $(\Delta(\cX),\Delta(\cX)^\perp)$ in $\cB\ltimes_\eta F$.

\begin{theorem}\label{thm:Delta-cotorsion}
If $(\cX,\cY)$ is a cotorsion pair in $\cB$, then
$
(\Delta(\cX),\Delta(\cX)^\perp)
$
is a cotorsion pair in $\cB\ltimes_\eta F$.
\end{theorem}

\begin{proof}
Since $\cX={}^\perp\cY$ and $\cY$ contains all injective objects, Lemma~\ref{lem:Delta-orthogonal} gives
$
\Delta(\cX)={}^\perp Z(\cY).
$
Therefore
\[
{}^\perp(\Delta(\cX)^\perp)
 ={}^\perp\bigl(({}^\perp Z(\cY))^\perp\bigr)
 ={}^\perp Z(\cY)
 =\Delta(\cX).
\]
Thus $(\Delta(\cX),\Delta(\cX)^\perp)$ is a cotorsion pair.
\end{proof}

\subsection{Heredity of the \texorpdfstring{$\Delta$}{Delta}-pair} In the case of the right trivial extension $\cB\ltimes F$, Hu \cite[Lemma 3.9]{Hu} proves that, under the condition $\mathbb{L}_1F(\cC)=0$,
$\Delta(\cC)$ is closed under the kernels of epimorphisms if and only if $\cC$ is closed under the kernels of epimorphisms. This implies heredity of
the cotorsion pair $(\Delta(\cX),\Delta(\cX)^\perp)$. For the $\eta$-extension $\cB\ltimes_\eta F$, this condition $\mathbb{L}_1F(\cC)=0$ is not sufficient. To address this case, we need the defect functor $\delta=C\Omega$. We shall show that the first left derived functor $\mathbb{L}_1\delta$ provides the precise heredity criterion for the general setting.

The only if part of \cite[Lemma 3.9]{Hu} holds also for the $\eta$-extension $\cB\ltimes_\eta F$. The proof is the same as that of \cite[Lemma 3.9]{Hu} and we omit it.

\begin{lemma}
Let $\cC$ be a class of objects in $\cB$. Assume that $\mathbb{L}_1F(\cC)=0$. If $\Delta(\cC)$ is closed under the kernels of epimorphisms, then so is $\cC$.
\end{lemma}

For the if part of \cite[Lemma 3.9]{Hu}, we replace the condition $\mathbb{L}_1F(\cC)=0$ with $\mathbb{L}_1\delta(\Delta(\cC))=0$.

\begin{lemma}\label{lem:ker-epi}
Let $\cC$ be a class of objects in $\cB$. Assume that $\mathbb{L}_1\delta(\Delta(\cC))=0$. If $\cC$ is closed under the kernels of epimorphisms, then so is $\Delta(\cC)$.
\end{lemma}
\begin{proof}
Let $0\rightarrow (X_1,f_1)\rightarrow (X_2,f_2)\rightarrow (X_3,f_3)\rightarrow 0$ be a short exact sequence with $(X_2,f_2), (X_3,f_3)\in \Delta(\cC)$.
We have a commutative diagram with exact rows
\[
\xymatrix{
0\ar[r] & \delta(X_1,f_1) \ar[r] \ar[d]^{\gamma_{(X_1,f_1)}} & \delta(X_2,f_2) \ar[r] \ar[d]^{\gamma_{(X_2,f_2)}} & \delta(X_3,f_3) \ar[r] \ar[d]^{\gamma_{(X_3,f_3)}} \ar[r] & 0\\
0\ar[r] & U(X_1,f_1) \ar[r] & U(X_2,f_2) \ar[r] & U(X_3,f_3) \ar[r]  & 0
}
\]
Here the top row is exact because $\mathbb{L}_1\delta(\Delta(\cC))=0$ and $(X_3,f_3)\in \Delta(\cC)$. Now applying the snake lemma, and noting that $\Coker(\gamma_{(X_i,f_i)})=C(X_i,f_i)$ for $i=1,2,3$ from \eqref{eq:delta-seq}, we have a long exact sequence
\[
\begin{aligned}
0\to\Ker\gamma_{(X_1,f_1)}&\to \Ker\gamma_{(X_2,f_2)} \to \Ker\gamma_{(X_3,f_3)} \\
&\to C(X_1,f_1) \to C(X_2,f_2) \to C(X_3,f_3) \to 0
\end{aligned}
\]
Since $(X_2,f_2), (X_3,f_3)\in \Delta(\cC)$, we know that $\gamma_{(X_2,f_2)}$ and $\gamma_{(X_3,f_3)}$ are monic and
$C(X_2,f_2),C(X_3,f_3)\in \cC$. The long exact sequence implies that
$\Ker\gamma_{(X_1,f_1)}=0$, i.e., $\gamma_{(X_1,f_1)}$ is monic,
and that the sequence $0\to C(X_1,f_1) \to C(X_2,f_2) \to C(X_3,f_3) \to 0$ is exact. Since $\cC$ is closed under the kernels of epimorphisms,
we have $C(X_1,f_1)\in \cC$.
Thus $(X_1,f_1)\in \Delta(\cC)$ by Lemma~\ref{lem:Com-ex}. This shows that $\Delta(\cC)$ is closed under  the kernels of epimorphisms.
\end{proof}

The following proposition shows that the two conditions $\mathbb{L}_1F(\cC)=0$ and $\mathbb{L}_1\delta(\Delta(\cC))=0$ are equivalent in the right trivial extension case.

\begin{proposition}
Let $\cC$ be a class of objects in $\cB$. If $\eta=0$, then $\mathbb{L}_1\delta(\Delta(\cC))=0$ if and only if $\mathbb{L}_1F(\cC)=0$.
\end{proposition}
\begin{proof}
Since $\eta=0$, for any object $(X,f)\in \cB\ltimes_\eta F$, we have
\[
\delta(X,f)=\Coker(\eta_X-F(f))\cong \Coker(F(f))\cong F(\Coker(f))=FC(X,f).
\]

Now assume that $\mathbb{L}_1F(\cC)=0$. Let
\begin{equation}\label{eq:ses}
0\rightarrow (X_1,f_1)\rightarrow (X_2,f_2)\rightarrow (X_3,f_3)\rightarrow 0
\end{equation}
be a short exact sequence
with $(X_3,f_3)\in \Delta(\cC)$. We have $\mathbb{L}_1C(X_3,f_3)=0$ by Lemma~\ref{lem:Com-ex}. Since $C(X_3,f_3)\in \cC$, we have $\mathbb{L}_1F(C(X_3,f_3))=0$. Applying the functors
$C$ and $F$ successively to the preceding short exact sequence \eqref{eq:ses}, by Lemma~\ref{lem:L_1-seq} we obtain an exact sequence
\[
0\rightarrow \delta(X_1,f_1)\rightarrow \delta(X_2,f_2)\rightarrow \delta(X_3,f_3)\rightarrow 0.
\]
This shows that $\delta$ preserves the exactness of short exact sequences with $(X_3,f_3)$ as the right end. Thus $\mathbb{L}_1\delta(X_3,f_3)=0$. Hence
$\mathbb{L}_1\delta(\Delta(\cC))=0$.

Conversely, assume that $\mathbb{L}_1\delta(\Delta(\cC))=0$. Let $0\to Y_1 \to Y_2 \to Y_3 \to 0$ be a short exact sequence with $Y_3\in \cC$. By Lemma~\ref{lem:E_Y-exact},
we obtain an exact sequence
\[
0\to Z(Y_1)\to E_{Y_2} \to T(Y_3) \to 0,
\]
and applying the functor $C$ to it recovers the original sequence.
Since $T(Y_3)\in \Delta(\cC)$ and $\mathbb{L}_1\delta(\Delta(\cC))=0$, we have $\mathbb{L}_1\delta(T(Y_3))=0$. It follows that the sequence
\[
0\to \delta(Z(Y_1))\to \delta(E_{Y_2}) \to \delta(T(Y_3))\to 0
\]
is exact. Since $\delta\cong FC$, the sequence $0\to F(Y_1)\to F(Y_2)\to F(Y_3)\to 0$ is exact. Thus $\mathbb{L}_1F(Y_3)=0$ by lemma~\ref{lem:L_1-seq}. Hence $\mathbb{L}_1F(\cC)=0.$
\end{proof}

We now give a criterion for the heredity of the cotorsion pair $(\Delta(\cX),\Delta(\cX)^\perp)$.

\begin{theorem}\label{thm:Delta-heredity}
Assume that $(\cX,\cY)$ is a hereditary cotorsion pair in $\cB$. Then the cotorsion pair $(\Delta(\cX),\Delta(\cX)^\perp)$ is hereditary if
and only if $\mathbb{L}_1\delta(\Delta(\cX))=0$.
\end{theorem}
\begin{proof}
Since $\cB\ltimes_\eta F$ has enough projective objects,
the if part is straightforward from Lemma~\ref{lem:ker-epi} and Proposition~\ref{prop:one-sided-criteria}~(1). Now consider the only if part. Assume that $(\Delta(\cX),\Delta(\cX)^\perp)$ is hereditary.
Thus $\Delta(\cX)$ is closed under the kernels of epimorphisms. For any $A=(X,f)\in \Delta(\cX)$,
since the $\eta$-extension $\cB\ltimes_\eta F$ has enough projective objects,
we can choose an exact sequence $0\rightarrow K\xrightarrow{i} P\rightarrow A\rightarrow 0$ with $P$ projective.
Note that $P\in \Delta(\cX)$ as the left half of a cotorsion pair contains all projective objects. Thus $K\in \Delta(\cX)$ and hence $\gamma_K$ is monic by lemma~\ref{lem:Com-ex}.

On the one hand, we have the following exact sequence
\[
0=\mathbb{L}_1\delta(P) \rightarrow \mathbb{L}_1\delta(A) \rightarrow \delta(K)\xrightarrow{\delta(i)} \delta(P) \rightarrow \delta(A) \rightarrow 0.
\]
On the other hand, we have a commutative diagram
\[
\xymatrix{
\delta(K)\ar[d]^{\gamma_K}\ar[r]^{\delta(i)} & \delta(P) \ar[d]^{\gamma_P}\\
U(K) \ar[r]^{U(i)} & U(P)
}
\]
Thus $\gamma_P\delta(i)=U(i)\gamma_K$, which implies that $\delta(i)$ is a monomorphism as $\gamma_K$ and $U(i)$ are monic. Combining both, we obtain that $\mathbb{L}_1\delta(A)=0$.
Therefore, $\mathbb{L}_1\delta(\Delta(\cX))=0$.
\end{proof}

\subsection{The \texorpdfstring{$U^{-1}$}{U-inverse}-pair and its relation to the \texorpdfstring{$\Delta$}{Delta}-pair} Following \cite{Hu}, for a class $\cC$ of objects in $\cB$, we define
\[
U^{-1}(\cC)=\{(X,f)\in \cB\ltimes_\eta F\mid X\in \cC\}.
\]
Under the condition $\mathbb{L}_1F(\cC)=0$, Hu \cite{Hu} gives a characterization of $U^{-1}(\cC^\perp)$ and proves that if $(\cX,\cY)$ is a cotorsion pair then
so is $({}^\perp U^{-1}(\cY),U^{-1}(\cY))$. The heredity criterion for $({}^\perp U^{-1}(\cY),U^{-1}(\cY))$ is also given. These results are obtained under the standing assumption that $\eta$ is associative and nilpotent. We observe, however, that the nilpotence is superfluous: all the proofs of these results in \cite{Hu} carry over to the associative case without modification. Accordingly, we record the results in this more general setting for future reference.

\begin{lemma}\label{lem:T-orthogonal}
{\rm (\cite[Lemma 3.6]{Hu})}
Let $\cC$ be a class of objects in $\cB$. If $\mathbb{L}_1F(\cC)=0$, then
$
T(\cC)^\perp=U^{-1}(\cC^\perp).
$
\end{lemma}

\begin{theorem}\label{thm:U-cotorsion}
{\rm (\cite[Theorem 3.5]{Hu})}
Let $\cB\ltimes_\eta F$ be an $\eta$-extension of $\cB$. Assume that $(\cX,\cY)$ is a cotorsion pair in $\cB$. If $\mathbb{L}_1F(\cX)=0$, then
$
({}^\perp U^{-1}(\cY),U^{-1}(\cY))
$
is a cotorsion pair in $\cB\ltimes_\eta F$. It is hereditary if and only if $(\cX,\cY)$ is hereditary.
\end{theorem}

The following proposition compares the two cotorsion pairs
$
(\Delta(\cX),\Delta(\cX)^\perp)
$ and
$
({}^\perp U^{-1}(\cY),U^{-1}(\cY)).
$

\begin{proposition}\label{prop:relation}
Suppose that $(\cX,\cY)$ is a cotorsion pair in $\cB$ and $\mathbb{L}_1F(\cX)=0$. Then
$
\Delta(\cX)^\perp\subseteq U^{-1}(\cY),
$
or equivalently,
$
{}^\perp U^{-1}(\cY)\subseteq\Delta(\cX).
$
\end{proposition}

\begin{proof}
Since $\mathbb{L}_1F(\cX)=0$, both $(\Delta(\cX),\Delta(\cX)^\perp)$ and $({}^\perp U^{-1}(\cY),U^{-1}(\cY))$ are cotorsion pairs. By Lemma~\ref{lem:T-orthogonal},
$
U^{-1}(\cX^\perp)=T(\cX)^\perp.
$
Since $T(\cX)\subseteq \Delta(\cX)$, taking right orthogonals yields
$
\Delta(\cX)^\perp\subseteq U^{-1}(\cY).
$
Then taking left orthogonals gives
$
{}^\perp U^{-1}(\cY)\subseteq\Delta(\cX).
$
\end{proof}

The proof of Proposition~\ref{prop:relation} shows that if $\Delta(\cX)=T(\cX)$, then the two cotorsion pairs coincide. Now we consider when the condition $\Delta(\cX)=T(\cX)$ holds.
Recall from Subsection~2.3 that for any $A\in \cB\ltimes_\eta F$ we have a short exact sequence $0\rightarrow \Omega(A)\xrightarrow{\omega_{A}} TU(A)\xrightarrow{\varepsilon_{A}} A\rightarrow 0$ and applying the functor $C$ to it gives an exact sequence
$\delta(A)\xrightarrow{\gamma_A} U(A)\xrightarrow{\rho_A} C(A)\rightarrow 0$.

\begin{lemma}\label{lem:T-equi}
Assume that $\eta$ is nilpotent. For $A=(X,f)\in \cB\ltimes_\eta F$, the following statements are equivalent:
\begin{enumerate}[label=\textup{(\arabic*)}]
\item $0\rightarrow \delta(A)\xrightarrow{\gamma_A} U(A)\xrightarrow{\rho_A} C(A)\rightarrow 0$ is split exact;
\item $A\cong TC(A)$;
\item there exists an object $Y\in \cB$ such that $A\cong T(Y)$.
\end{enumerate}
\end{lemma}
\begin{proof}
{\rm (1) $\Rightarrow$ (2):} Since the sequence is split exact, choose a morphism $s_A:C(A)\rightarrow U(A)$ such that $\rho_As_A=1$. Consider the morphism
$\alpha=\varepsilon_AT(s_A)\colon TC(A)\xrightarrow{T(s_A)}TU(A)\xrightarrow{\varepsilon_A}A$. We claim that $\alpha$ is an isomorphism. Indeed, noting that $CT=\Id_{\cB}$, we have
$C(\alpha)=C(\varepsilon_A)CT(s_A)=\rho_As_A=1$. By \cite[Lemma~2.3]{Marmaridis} or \cite[Lemma 5.3]{Beligiannis}, nilpotence of $\eta$ implies that $\Ker C=\{(X_1,f_1)\in \cB\ltimes_\eta F\mid
C(X_1,f_1)=0\}=0$. Since $\gamma_A$ is monic, by Lemma~\ref{lem:Com-ex} we have $\mathbb{L}_1C(A)=0$. Thus by \cite[Corollary~5.10 (i),(iii)]{Beligiannis} we obtain that $\alpha$ is an isomorphism.
Therefore, $A\cong TC(A)$.

{\rm (2) $\Rightarrow$ (3):} It is immediate.

{\rm (3) $\Rightarrow$ (1):} It suffices to show that $0\rightarrow \delta(T(Y))\xrightarrow{\gamma_{T(Y)}} U(T(Y))\xrightarrow{\rho_{T(Y)}} C(T(Y))\rightarrow 0$
is split exact. First, since $\Com(T(Y))$ is exact, $\gamma_{T(Y)}$ is monic. Second,
the morphism $s=\binom{1}{0}\colon Y\to Y\oplus F(Y)$ satisfies $\rho_{T(Y)}s=1$, where $\rho_{T(Y)}=(1,0)\colon Y\oplus F(Y)\to Y$. Therefore, the sequence is split exact.
\end{proof}
\begin{remark}
{\rm The implication (3) $\Rightarrow$ (2) is clear: if $A\cong T(Y)$, then applying the functor $C$ we obtain $C(A)\cong CT(Y)=Y$. Thus the equivalence of (2) and (3) follows naturally.}
\end{remark}

\begin{proposition}\label{prop:Delta-induced}
Let $\cC$ be a class of objects in $\cB$. Assume that $\eta$ is nilpotent. Then $\Delta(\cC)=T(\cC)$ if and only if the canonical sequence
\begin{equation}\label{eq:T-canonical seq}
0\longrightarrow \delta(A)\xlongrightarrow{\gamma_A}U(A)\xlongrightarrow{\rho_A}C(A)\longrightarrow0
\end{equation}
is split exact for every $A\in\Delta(\cC)$. In particular, the equality holds if
$\Ext^1_\cB(C(A),\delta(A))=0$ for any $A\in\Delta(\cC)$.
\end{proposition}
\begin{proof}
Assume first that $\Delta(\cC)=T(\cC)$ and let $A\in\Delta(\cC)$. Then $A\cong T(Y)$ for some $Y\in\cC$. Lemma~\ref{lem:T-equi} shows that the canonical sequence \eqref{eq:T-canonical seq} of $A$ is split exact.

Conversely, assume that the canonical sequence \eqref{eq:T-canonical seq} is split exact for every $A\in\Delta(\cC)$. By Lemma~\ref{lem:T-equi}, $A\cong T(C(A))$; since $C(A)\in\cC$, this gives $A\in T(\cC)$. Hence $\Delta(\cC)\subseteq T(\cC)$. The reverse inclusion always holds, so equality follows.

Finally, for $A\in \Delta(\cC)$, by lemma~\ref{lem:Com-ex}, $\gamma_A$ is monic and hence the sequence \eqref{eq:T-canonical seq} is exact. $\Ext^1_\cB(C(A),\delta(A))=0$ implies that the exact
sequence \eqref{eq:T-canonical seq} is split.
\end{proof}

\begin{theorem}\label{thm:two-equal}
Let $(\cX,\cY)$ be a cotorsion pair in $\cB$. Assume that $\eta$ is nilpotent and $\mathbb{L}_1F(\cX)=0$. If $\Ext^1_{\cB}(C(A),\delta(A))=0$ for any $A\in \Delta(\cX)$, then
$(\Delta(\cX),\Delta(\cX)^\perp)=({}^\perp U^{-1}(\cY),U^{-1}(\cY))=(T(\cX),U^{-1}(\cY))$.
\end{theorem}
\begin{proof}
Proposition~\ref{prop:Delta-induced} gives $\Delta(\cX)=T(\cX)$. By Lemma~\ref{lem:T-orthogonal},
$
T(\cX)^\perp=U^{-1}(\cX^\perp)=U^{-1}(\cY).
$
Therefore
$
\Delta(\cX)^\perp=T(\cX)^\perp=U^{-1}(\cY).
$
Since $(\Delta(\cX),\Delta(\cX)^\perp)$ and $({}^\perp U^{-1}(\cY),U^{-1}(\cY))$ are cotorsion pairs, $\Delta(\cX)={}^\perp(\Delta(\cX)^\perp)={}^\perp U^{-1}(\cY)$,
and the asserted equalities follow.
\end{proof}

\begin{remark}\label{rem:equal-reduce}
{\rm
When $\eta=0$, for any $A=(X,f)\in \Delta(\cX)$, since $F$ is right exact, it follows that $\delta(A)=\Coker(\eta_{X}-F(f))\cong \Coker(F(f))\cong F(\Coker(f))=FC(A)$. Consequently,
the condition $\Ext^1_{\cB}(C(A),\delta(A))=0$ for any $A\in \Delta(\cX)$ follows from the condition $\Ext^1_{\cB}(X,F(X))=0$ for any $X\in \cX$. Conversely,
assume that $\Ext^1_{\cB}(C(A),\delta(A))=0$ for any $A\in \Delta(\cX)$. For any $X\in \cX$, we have $T(X)\in \Delta(\cX)$, and hence
$\Ext^1_{\cB}(X,F(X))\cong \Ext^1_{\cB}(C(T(X)),\delta(T(X)))=0$. Therefore the two conditions are equivalent in this case and Theorem~\ref{thm:two-equal} reduces to that of \cite[Theorem~3.11]{Hu}.
}
\end{remark}

\subsection{Completeness of cotorsion pairs}

We shall establish completeness criteria of the two cotorsion pairs $(\Delta(\cX),\Delta(\cX)^\perp)$ and $({}^\perp U^{-1}(\cY),U^{-1}(\cY))$. In \cite{Hu}, only the completeness
of the latter is considered in the right trivial extension case; see \cite[Theorem~3.13]{Hu}.

For ease of exposition, we introduce the following notion.

\begin{definition}
{\rm (\cite[Definition 3.1.1]{GobelTrlifaj})}
Let $\cA$ be an abelian category, $\cC$ a class of objects in $\cA$ and $X$ an object of $\cA$. $X$ is said to be finitely $\cC$-filtered if
there is a finite sequence of monomorphisms
\[
0=X_0\rightarrowtail X_1\rightarrowtail \cdots \rightarrowtail X_n=X
\]
such that the quotient $X_{i+1}/X_i:=\Coker(X_i\rightarrowtail X_{i+1})\in \cC$ for every $i$.
\end{definition}

Recall that $\radf$ is the radical functor of Beligiannis in Subsection~2.3. The following two lemmas are implicit in \cite{Beligiannis}. See also \cite{LG}.

\begin{lemma}\label{lem:nilp-radf}
If $\eta$ is $N$-nilpotent, then $\radf^{N+1}=0$.
\end{lemma}
\begin{proof}
Let $(X,f)\in\cB\ltimes_\eta F$. Since $\eta$ is $N$-nilpotent, we have
\[
\eta_X\eta_{F(X)}\cdots\eta_{F^{N-1}(X)}=0.
\]
Associativity of $\eta$ and the relation $fF(f)=f\eta_X$ imply
$
fF(f)\cdots F^N(f)=0.
$
Set $X_0=X$ and $f_0=f$. For $s\geq 0$, let
$
F(X_s)\xrightarrow{k_s}\operatorname{Im}(f_s)\xrightarrow{i_s}X_s
$
be the canonical factorization of $f_s$ and set
$
X_{s+1}=\Ima(f_s),
$
$
f_{s+1}=k_sF(i_s).
$
By the definition of $\radf$, $\radf^s(X,f)=(X_s,f_s)$ for any $s\geq 1$.

Moreover, a direct computation yields
\[
fF(f)\cdots F^N(f)
 =i_0i_1\cdots i_{N-1}
 f_NF(k_{N-1})\cdots F^{N-1}(k_1)F^N(k_0).
\]
Each $i_s$ is monic, and each $F^{N-s}(k_s)$ is epic because $F$ is right exact. Thus $f_N=0$, so $(X_N,f_N)=Z(X_N)$ and $X_{N+1}=\Ima(0)=0$.
Hence $\radf^{N+1}(X,f)=0$. This shows that $\radf^{N+1}=0$.
\end{proof}

\begin{lemma}\label{lem:nilp-Z-filt}
If $\eta$ is nilpotent, then every object $(X,f)\in \cB\ltimes_\eta F$ is finitely $Z(\cB)$-filtered.
\end{lemma}
\begin{proof}
Assume that $\eta$ is $N$-nilpotent. Then by Lemma~\ref{lem:nilp-radf}, $\radf^{N+1}=0$. For every object $(X,f)\in \cB\ltimes_\eta F$,
and $0\leq s\leq N$, by \eqref{eq:radical-seq} we have an exact sequence
\[
0\to \radf^{s+1}(X,f)\to \radf^s(X,f)\to ZC(\radf^s(X,f))\to 0.
\]
Since $\radf^{N+1}=0$, we obtain a finite sequence of monomorphisms
\[
0=\radf^{N+1}(X,f)\rightarrowtail \radf^{N}(X,f)\rightarrowtail \cdots \rightarrowtail \radf^{0}(X,f)=(X,f),
\]
with the quotient $\radf^{s}(X,f)/\radf^{s+1}(X,f)\cong ZC(\radf^s(X,f))\in Z(\cB)$ for every $s$. Therefore, $(X,f)$ is  finitely $Z(\cB)$-filtered.
\end{proof}

The following lemma is crucial, as it allows us to obtain sufficient conditions for the completeness of cotorsion pairs.

\begin{lemma}\label{lem:nilp-Z(Y)}
Let $\cY$ be a class in $\cB$ and let $Y\in \cY$. For an exact sequence $0\rightarrow B\xrightarrow{j} Y \xrightarrow{p} X\rightarrow 0$, assume that $0\rightarrow Z(B)\rightarrow E_Y\rightarrow T(X)\rightarrow 0$ is the sequence constructed as in Lemma~\ref{lem:E_Y-exact}. If $\eta$ is $N$-nilpotent and $C\radf^nT(X)\in \cY$ for
every $1\leq n\leq N$, then $E_Y$ is finitely $Z(\cY)$-filtered.
\end{lemma}
\begin{proof}
By Lemma~\ref{lem:nilp-radf}, $\radf^{N+1}=0$. A direct computation gives that
$\radf(E_Y)\cong (F(X),\eta_X)\cong \radf T(X)$. Applying the functor $\radf^{n-1}$ to this isomorphism, we obtain $\radf^n(E_Y)\cong \radf^n T(X)$ for every $n\geq 1$.
For every $0\leq n\leq N$, by \eqref{eq:radical-seq} we have the following exact sequence
\[
0\to \radf^{n+1}(E_Y)\to \radf^n(E_Y) \to ZC\radf^n(E_Y)\to 0.
\]
Note that $C(E_Y)=Y\in \cY$ and $C\radf^n(E_Y)\cong C\radf^nT(X)\in \cY$ for
every $1\leq n\leq N$. Thus we obtain a finite sequence of monomorphisms
\[
0=\radf^{N+1}(E_Y)\rightarrowtail \radf^N(E_Y)\rightarrowtail \cdots \rightarrowtail \radf^0(E_Y)=E_Y,
\]
with all quotients $ZC\radf^n(E_Y)\in Z(\cY)$. This shows that $E_Y$ is finitely $Z(\cY)$-filtered.
\end{proof}

Noting that $\cB\ltimes_\eta F$ has enough projective objects, we have the following completeness criterion for $(\Delta(\cX),\Delta(\cX)^\perp)$.

\begin{theorem}\label{thm:completeness}
Let $(\cX,\cY)$ be a hereditary complete cotorsion pair in $\cB$. Assume that $\eta$ is $N$-nilpotent. If $\mathbb{L}_1\delta(\Delta(\cX))=0$ and $C\radf^nT(\cX)\subseteq \cY$ for all $1\leq n\leq N$, then $(\Delta(\cX),\Delta(\cX)^\perp)$ is
a hereditary complete cotorsion pair in $\cB\ltimes_\eta F$.
\end{theorem}
\begin{proof}
By Theorems~\ref{thm:Delta-cotorsion} and \ref{thm:Delta-heredity}, $(\Delta(\cX),\Delta(\cX)^\perp)$ is a hereditary cotorsion pair.
Now consider the completeness.

Let $B\in \cB$. Since $(\cX, \cY)$ is a complete cotorsion pair, there is an exact sequence $0\to B\to Y\to X\to 0$ with $Y\in \cY, X\in \cX$. By Lemma~\ref{lem:E_Y-exact}, we
have an exact sequence
\[
0\to Z(B)\to E_Y \to T(X)\to 0,
\]
where the underlying object of $E_Y$ in $\cB$ is $Y\oplus F(X)$.
By Lemma~\ref{lem:Delta-orthogonal}, $Z(\cY)\subseteq ({}^\perp Z(\cY))^\perp = \Delta(\cX)^\perp$.

Since $\eta$ is $N$-nilpotent and $C\radf^nT(\cX)\subseteq \cY$ for all $1\leq n\leq N$,
by Lemma~\ref{lem:nilp-Z(Y)} $E_Y$ is finitely $Z(\cY)$-filtered. Note that the right half $\Delta(\cX)^\perp$ of the cotorsion
pair $(\Delta(\cX),\Delta(\cX)^\perp)$ is closed under extensions. Thus $E_Y\in \Delta(\cX)^\perp$. On the other hand, $T(X)\in \Delta(\cX)$. Hence
$Z(B)$  has a special $\Delta(\cX)^{\perp}$-preenvelope.

From Lemma~\ref{lem:nilp-Z-filt}, we know that every object $(X,f)\in \cB\ltimes_\eta F$ is
finitely $Z(\cB)$-filtered. Thus $(X,f)$ has a special $\Delta(\cX)^\perp$-preenvelope by Proposition~\ref{prop:extension-approx} (2), and hence
the cotorsion pair $(\Delta(\cX),\Delta(\cX)^\perp)$ is complete by Proposition~\ref{prop:one-sided-complete} (1).
\end{proof}

The following completeness criterion for $({}^\perp U^{-1}(\cY), U^{-1}(\cY))$ is contained in \cite{LG}.

\begin{theorem}\label{thm:pre-U-completeness}
{\rm (See \cite[Theorem~4.13]{LG})}
Let $(\cX,\cY)$ be a hereditary complete cotorsion pair in $\cB$. Assume that $\eta$ is nilpotent. If $\mathbb{L}_1F(\cX)=0$ and $F(\cX)\subseteq \cY$, then $({}^\perp U^{-1}(\cY), U^{-1}(\cY))$ is
a hereditary complete cotorsion pair in $\cB\ltimes_\eta F$.
\end{theorem}

In general, the fact that $\cB$ has enough injective objects does not guarantee that the extension category $\cB\ltimes_\eta F$ does. However, if $\cB\ltimes_\eta F$
does have enough injective objects, we can deduce the completeness of $({}^\perp U^{-1}(\cY), U^{-1}(\cY))$ without the nilpotence assumption, as follows.

\begin{theorem}\label{thm:U-completeness}
Let $(\cX,\cY)$ be a complete cotorsion pair in $\cB$. Assume that $\cB\ltimes_\eta F$ has enough injective objects. If $\mathbb{L}_1F(\cX)=0$ and $F(\cX)\subseteq \cY$, then $({}^\perp U^{-1}(\cY), U^{-1}(\cY))$ is
a complete cotorsion pair in $\cB\ltimes_\eta F$.
\end{theorem}
\begin{proof}
By Theorem~\ref{thm:U-cotorsion}, $({}^\perp U^{-1}(\cY), U^{-1}(\cY))$ is a cotorsion pair. It remains to prove completeness.
Let $(B,f)\in \cB\ltimes_\eta F$. Since $(\cX,\cY)$ is a complete cotorsion pair, there is a short exact sequence $0\to Y\to X\xrightarrow{p} B\to 0$ with $X\in \cX$ and $Y\in \cY$.
$(p, f\circ F(p))\colon X\oplus F(X)\to B$ gives a morphism from $T(X)$ to $(B,f)$ because the following diagram commutes
\[
\xymatrix@C=2.5pc{
F(X)\oplus F^2(X)\ar[rr]^{\ \ \ \ \ (F(p), F(f)\circ F^2(p))}\ar[d]^{t_X=\big(\begin{smallmatrix}0 & 0\\ 1 & \eta_X\end{smallmatrix}\big)} && F(B)\ar[d]^{f}\\
X\oplus F(X) \ar[rr]^{(p, f\circ F(p))} && B
}
\]
Since $p$ is an epimorphism, it's easy to see that $(p, f\circ F(p))$ is also an epimorphism. Thus we get a short exact sequence
\begin{equation}\label{eq:special-precover}
0\to K\to T(X)\xrightarrow{(p, f\circ F(p))}(B,f)\to 0.
\end{equation}
Applying the exact functor $U$, we have an exact sequence
$
0\to U(K)\to X\oplus F(X)\xrightarrow{(p, f\circ F(p))}B\to 0.
$
Consider the following commutative diagram
\[
\xymatrix{
0 \ar[r] & X \ar[r]^{\binom{1}{0}\ \ \ \ \ \ }\ar[d]^{p} & X\oplus F(X) \ar[r]^{\ \ \ (0,1)}\ar[d]^{(p, f\circ F(p))} & F(X) \ar[r]\ar[d] & 0\\
0 \ar[r] & B \ar[r]^{1} & B \ar[r] & 0 \ar[r] & 0
}
\]
By the snake lemma, we have a short exact sequence $0\to Y\to U(K)\to F(X)\to 0$. $F(\cX)\subseteq \cY$ implies that $F(X)\in \cY$. Since $\cY$ is closed under extensions,
$U(K)\in \cY$. So $K\in U^{-1}(\cY)$. By Lemma~\ref{lem:T-orthogonal}, $T(X)\in {}^\perp(T(\cX)^\perp)={}^\perp U^{-1}(\cY)$. The sequence \eqref{eq:special-precover}
shows that $(B,f)$ has a special ${}^\perp U^{-1}(\cY)$-precover. Thus ${}^\perp U^{-1}(\cY)$ is special precovering.
Since $\cB\ltimes_\eta F$ has enough injective objects, it follows from Proposition~\ref{prop:one-sided-complete} (2) that $({}^\perp U^{-1}(\cY), U^{-1}(\cY))$ is
complete.
\end{proof}

\subsection{Cotorsion pairs in \texorpdfstring{$\zeta$}{zeta}-coextensions}
Let $\cB$ be an abelian category with enough projective and injective objects, $G\colon\cB\to\cB$ a left exact functor, and $\zeta\colon G\to G^2$ a natural transformation. The category $G\rtimes_\zeta\cB$, called a $\zeta$-coextension of $\cB$, is defined as follows. The objects of $G\rtimes_\zeta\cB$ are pairs $[X,h]$, where $X$ is an object of $\cB$ and $h\colon X\to G(X)$ is a morphism of $\cB$ satisfying
$
G(h)\circ h=\zeta_X\circ h.
$
A morphism from $[X,h]$ to $[Y,k]$ is a morphism $\alpha\colon X\to Y$ of $\cB$ such that
$
k\circ \alpha=G(\alpha)\circ h.
$

From \cite{Marmaridis}, the $\zeta$-coextension $G\rtimes_\zeta\cB$ is an abelian category. The sequence
$[X_1,h_1]\xrightarrow{\alpha} [X_2,h_2]\xrightarrow{\beta} [X_3,h_3]$ in $G\rtimes_\zeta\cB$ is exact if and only if the underlying sequence
$X_1\xrightarrow{\alpha} X_2\xrightarrow{\beta} X_3$ in $\cB$ is exact. When $\zeta=0$, we write $G\rtimes\cB$ and call it the left trivial extension of $\cB$ by $G$ \cite{FossumGriffithReiten}.

The natural transformation $\zeta\colon G\to G^2$ is \emph{coassociative} if
$
G(\zeta_X)\circ \zeta_X=\zeta_{G(X)}\circ\zeta_X
$
for every object $X\in\cB$. It is \emph{nilpotent} if there exists $N\ge1$ such that
$
\zeta_{G^{N-1}(X)}\circ\cdots\circ\zeta_{G(X)}\circ\zeta_{X}=0,
$
for every object $X\in\cB$. In this case, we say that $\zeta$ is \emph{$N$-nilpotent}. Throughout this subsection, $\zeta$ is assumed to be coassociative.

From \cite{Beligiannis}, there are the following standard functors:

\[
\xymatrix@C=3pc{
\cB \ar@<1ex>[r]^{Z\ \ \ } & G\rtimes_\zeta\cB \ar@<1ex>[l]^{K\ \ \ } \ar@<1ex>[r]^{\ \ \ U}
& \cB \ar@<1ex>[l]^{\ \ \ H}
}
\]

\begin{itemize}
 	\item $U[X,f]=X$ and $U(\alpha)=\alpha$;
 	\item $Z(X)=[X,0]$ and $Z(\alpha)=\alpha$;
 	\item $H(X)=\left(G(X)\oplus X, s_X =\left(\begin{smallmatrix}
 		\zeta_X &0\\
 		1&0
 	\end{smallmatrix}\right)\right)$ and $H(\alpha)=\begin{pmatrix}
 	G(\alpha)&0\\
 	0& \alpha
 	\end{pmatrix}$;
 	\item $K[X,f]=\Ker(f)$ and $K(\alpha)$ is the morphism induced by the kernel.
 \end{itemize}
Note that the coassociativity of $\zeta$ guarantees that
$H(X)$ belongs to $G\rtimes_\zeta \cB$.

The pairs $(U,H)$ and $(Z,K)$ are adjoint pairs of functors. We have the identifications: $KZ=\mathrm{Id}_{\cB}$, and $KH=\mathrm{Id}_{\cB}$.
The functors $U$ and $Z$ are exact.
If $\cB$ has enough injective objects, then so does $G\rtimes_\zeta \cB$.

For $A=[X,h]$ we write
\[
\CoCom(A):
X\xrightarrow{h}G(X)
\xrightarrow{\,\zeta_X-G(h)\,}G^2(X).
\]
For a class $\cC\subseteq\cB$, define
\[
\nabla(\cC)=
\{A=[X,h]\mid \CoCom(A)\text{ is exact and }K(A)\in\cC\}.
\]

The unit of $(U,H)$ gives a natural short exact sequence
\[
0\rightarrow [X,h]\xrightarrow{\binom{h}{1}} HU[X,h] \xrightarrow{(1,-h)} D[X,h]\rightarrow 0,
\]
where $D[X,h]=[G(X),\zeta_X-G(h)]$. We define the \emph{codefect functor} $\delta^\vee=KD\colon G\rtimes_\zeta \cB\to \cB$.
$\delta^\vee[X,h]=\Ker(\zeta_X-G(h))$.

The counit of $(Z,K)$ gives a natural short exact sequence
\[
0\rightarrow ZK[X,h]\xrightarrow{\kappa_{[X,h]}} [X,h] \rightarrow \Coker(\kappa_{[X,h]})\rightarrow 0.
\]
We define a functor $\coradf$ by $\coradf[X,h]= \Coker(\kappa_{[X,h]})$, called the \emph{coradical functor} of Beligiannis.

The following results are the dual versions of Theorems~\ref{thm:Delta-cotorsion}, \ref{thm:Delta-heredity}, \ref{thm:U-cotorsion}, \ref{thm:two-equal}, \ref{thm:completeness} and \ref{thm:U-completeness}.

\begin{theorem}\label{thm:nabla-cotorsion}
If $(\cX,\cY)$ is a cotorsion pair in $\cB$, then
$
({}^\perp\nabla(\cY),\nabla(\cY))
$
is a cotorsion pair in $G\rtimes_\zeta \cB$.
\end{theorem}

\begin{theorem}\label{thm:nabla-heredity}
Assume that $(\cX,\cY)$ is a hereditary cotorsion pair in $\cB$. Then the cotorsion pair $({}^\perp\nabla(\cY),\nabla(\cY))$ is hereditary if
and only if $\mathbb{R}^1\delta^\vee(\nabla(\cY))=0$.
\end{theorem}

\begin{theorem}\label{thm:coU-cotorsion}
{\rm (\cite[Theorem 3.21]{Hu})}
Let $(\cX,\cY)$ be a cotorsion pair in $\cB$. If $\mathbb{R}^1G(\cY)=0$, then
$
(U^{-1}(\cX),U^{-1}(\cX)^\perp)
$
is a cotorsion pair in $G\rtimes_\zeta \cB$. It is hereditary if and only if $(\cX,\cY)$ is hereditary.
\end{theorem}

\begin{theorem}\label{thm:cotwo-equal}
Let $(\cX,\cY)$ be a cotorsion pair in $\cB$. Assume that $\zeta$ is nilpotent and $\mathbb{R}^1G(\cY)=0$. If $\Ext^1_{\cB}(\delta^\vee(A),K(A))=0$ for any $A\in \nabla(\cY)$, then
$({}^\perp\nabla(\cY),\nabla(\cY))=(U^{-1}(\cX),U^{-1}(\cX)^\perp)=(U^{-1}(\cX),H(\cY))$.
\end{theorem}

\begin{theorem}\label{thm:cocompleteness}
Let $(\cX,\cY)$ be a hereditary complete cotorsion pair in $\cB$. Assume that $\zeta$ is $N$-nilpotent. If $\mathbb{R}^1\delta^\vee(\nabla(\cY))=0$ and $K\coradf^nH(\cY)\subseteq \cX$ for all $1\leq n\leq N$, then $({}^\perp\nabla(\cY),\nabla(\cY))$ is
a hereditary complete cotorsion pair in $G\rtimes_\zeta \cB$.
\end{theorem}

\begin{theorem}\label{thm:U-cocompleteness}
Let $(\cX,\cY)$ be a complete cotorsion pair in $\cB$. Assume that $G\rtimes_\zeta \cB$ has enough projective objects. If $\mathbb{R}^1G(\cY)=0$ and $G(\cY)\subseteq \cX$, then $(U^{-1}(\cX),U^{-1}(\cX)^\perp)$ is
a complete cotorsion pair in $G\rtimes_\zeta \cB$.
\end{theorem}

\subsection{Four cotorsion pairs in one category} Let $\cB\ltimes_\eta F$ be an $\eta$-extension. Assume that $F$ has a right adjoint functor $G$.
Then we can find a natural transformation $\zeta\colon G\to G^2$ such that $\cB\ltimes_\eta F\cong G\rtimes_\zeta \cB$ \cite[Proposition 4.2]{Marmaridis}.
Specifically, for any
two objects $X,Y\in \cB$, let $\phi_{X,Y}$ be the isomorphism from $\Hom_{\cB}(FX,Y)$ to $\Hom_{\cB}(X,GY)$, determined by the adjoint pair $(F,G)$, and let
$\psi_{X,Y}=\phi_{X,Y}^{-1}$. For notational simplicity, we omit the subscripts of $\phi$ and $\psi$, as this will cause no ambiguity. Define
$\zeta_X=G^2(\psi(1_{G(X)}))\circ \phi^2(\eta_{G(X)})\colon G(X)\to G^2(X)$.

\begin{proposition}
{\rm (See \cite[Section 4]{Marmaridis})}
For the above constructed $\zeta$, the following statements hold.
\begin{enumerate}[label=\textup{(\arabic*)}]
\item $\zeta$ is a natural transformation.
\item If $\eta$ is associative, then $\zeta$ is coassociative.
\item If $\eta$ is $N$-nilpotent, then $\zeta$ is $N$-nilpotent.
\item The assignment $(X,f)\mapsto [X,f^{\sharp}]$, where $f^{\sharp}=\phi(f)\colon X\to G(X)$, gives an isomorphism $\cB\ltimes_\eta F\cong G\rtimes_\zeta \cB$.
\end{enumerate}
\end{proposition}

Via this isomorphism, we identify $\cB\ltimes_\eta F$ with $G\rtimes_\zeta \cB$.

\begin{theorem}
Let $(\cX,\cY)$ be a cotorsion pair in $\cB$, $(F,G)$ an adjoint pair of endofunctors of $\cB$, and $\eta\colon F^2\to F$ an associative natural transformation.
Then 
\begin{enumerate}[label=\textup{(\arabic*)}]
\item $(\Delta(\cX), \Delta(\cX)^\perp)$ is a cotorsion pair in $\cB\ltimes_\eta F$.
\item $({}^\perp\nabla(\cY), \nabla(\cY))$ is a cotorsion pair in $\cB\ltimes_\eta F$.
\item If $\mathbb{L}_1F(\cX)=0$, then $({}^\perp U^{-1}(\cY), U^{-1}(\cY))$ is a cotorsion pair in $\cB\ltimes_\eta F$.
\item If $\mathbb{R}^1G(\cY)=0$, then $(U^{-1}(\cX),U^{-1}(\cX)^{\perp})$ is a cotorsion pair in $\cB\ltimes_\eta F$.
\end{enumerate}
\end{theorem}
\begin{proof}
Parts (1) and (3) are Theorems~\ref{thm:Delta-cotorsion} and \ref{thm:U-cotorsion}. Under the extension--coextension identification, parts (2) and (4) are Theorems~\ref{thm:nabla-cotorsion} and \ref{thm:coU-cotorsion}.
\end{proof}

Their heredity and completeness have been studied in the preceding subsections.

\section{Cotorsion pairs over \texorpdfstring{$\theta$}{theta}-extension rings}\label{section:theta}

In this section, we apply the general theory of Section~\ref{section:extension} to $\theta$-extension rings and obtain cotorsion pairs over $\theta$-extension rings.
Whenever there is no ambiguity, we freely use the notation from Sections~\ref{section:pre} and~\ref{section:extension}.

Let $R$ be a ring, and $L$ an $R$-bimodule. Let
$
\theta\colon L\otimes_RL\to L
$
be an associative $R$-bimodule homomorphism, that is, $\theta$ is an $R$-bimodule homomorphism such that
$
\theta\circ (\theta\otimes 1_L)=\theta\circ(1_L\otimes\theta).
$
The $\theta$-extension of $R$ by $L$, denoted by $R\ltimes_\theta L$, is defined to be the ring with underlying group $R\oplus L$
and the multiplication given by
\[
(r,l)(r',l')=(rr',\,rl'+lr'+\theta(l\otimes l')).
\]

Let $\Lambda=R\ltimes_\theta L$. We identify $L$ with the ideal $0\oplus L$ of $\Lambda$. Under this identification, $L$ has the induced left and right $\Lambda$-module structures.
In this section, whenever we regard $L$ as a left (or right)
$\Lambda$-module, we mean the module structure described above, and not the one induced by the ring homomorphism $\Lambda\to R: (r,l)\mapsto r$.
We also identify $R$ with the subring $R\oplus0$ of $\Lambda$.

A left $\Lambda$-module is equivalently an $R$-module $X$ together with an $R$-module homomorphism
$
f\colon L\otimes_RX\longrightarrow X
$
satisfying
\begin{equation}\label{eq:theta-module-assoc}
f\circ (1_L\otimes f)=f\circ(\theta\otimes1_X).
\end{equation}
Indeed, $f(l\otimes x)$ is the action of $(0,l)$ on $x$, and \eqref{eq:theta-module-assoc} is exactly associativity of this action. Consequently,
we have an isomorphism of categories
\[
\Lambda\text{-Mod}\cong (R\text{-Mod})\ltimes_\eta F,
\]
where $F=L\otimes_R-$ and $\eta$ is defined by $\eta_X=\theta\otimes1_X.$ Since $\theta$ is associative, $\eta$ is associative.
Under this isomorphism the forgetful functor $U$ is restriction of scalars along $R\hookrightarrow\Lambda$ and the induction functor is $T=\Lambda\otimes_R-$.

The functor $F=L\otimes_R-$ has a right adjoint
$
G=\Hom_R(L,-).
$
As shown in Subsection~3.6, there exists a natural transformation $\zeta\colon G\to G^2$ satisfying $(R\text{-Mod})\ltimes_\eta F\cong G\rtimes_\zeta (R\text{-Mod})$.
Indeed, $\zeta_X\colon \Hom_R(L,X)\to \Hom_R(L,\Hom_R(L,X))$ is just $\zeta_X(g)=g\circ \theta$ under the natural identification
$
\Hom_R(L,\Hom_R(L,X))\cong \Hom_R(L\otimes_R L, X).
$
Thus we have
\[
\Lambda\Mod\cong G\rtimes_\zeta (R\text{-Mod}).
\]
Under this isomorphism, the coinduction functor is $H=\Hom_R(\Lambda,-)$.

For a $\Lambda$-module $A=(X,f)$, where $X$ is the underlying $R$-module and $f\colon L\otimes_R X\to X$ is the structure map, define an $R$-module homomorphism
\[
f^\sharp\colon X\longrightarrow\Hom_R(L,X),
\]
by $f^\sharp(x)(l)=f(l\otimes x)=lx.$
The object in $G\rtimes_\zeta R\text{-Mod}$ corresponding to the $\Lambda$-module $A$ is $[X,f^\sharp]$. Using the isomorphisms
\[
\Lambda\text{-Mod}\cong (R\text{-Mod})\ltimes_\eta F\cong G\rtimes_\zeta (R\text{-Mod}),
\]
we identify these three categories.

\begin{proposition}\label{prop:theta-functors}
Let $A$ be a $\Lambda$-module. Under the identifications above, the functors of  Sections~{\ref{section:pre}} and~{\ref{section:extension}} have the following concrete descriptions:
\begin{align*}
C(A)&= A/LA, & \delta(A)&\cong L\otimes_\Lambda A,\\
K(A)&=\Ann_A(L), & \delta^\vee(A)&\cong\Hom_\Lambda(L,A).
\end{align*}
Here, $\Ann_A(L)=\{a\in A\mid la=0\text{ for all } l\in L\}$. Moreover, the morphism
$
\gamma_A\colon\delta(A)\to U(A)
$
corresponds to the multiplication map
\begin{equation}\label{eq:theta-mu}
\mu_A\colon L\otimes_\Lambda A\longrightarrow A,
\qquad l\otimes a\longmapsto la,
\end{equation}
and the morphism $\gamma^\vee_A\colon U(A)\to \delta^\vee(A)$ corresponds to
\begin{equation}\label{eq:theta-nu}
\mu^\vee_A\colon A\longrightarrow\Hom_\Lambda(L,A),
\qquad \mu^\vee_A(a)(l)=la.
\end{equation}
\end{proposition}

\begin{proof}
Write $A=(X,f)$. Since $LA=\Ima(f)$, we have
\[
C(A)=\Coker(f)= A/LA.
\]
Likewise, $K(A)=\Ker(f^\sharp)$ consists precisely of the elements annihilated by every $l\in L$; thus $K(A)=\Ann_A(L)$.

By the definition of the defect functor,
\[
\delta(A)=\Coker\bigl(\theta\otimes1_X-1_L\otimes f\bigr).
\]
On the other hand, $L\otimes_\Lambda A$ can be obtained from $L\otimes_RX$ by imposing the balancing relations
\[
\theta(l\otimes l')\otimes x
=l\otimes f(l'\otimes x)\qquad \forall l,l^\prime\in L, x\in X.
\]
Thus $\delta(A)\cong L\otimes_\Lambda A$, and the morphism $\gamma_A$ is exactly \eqref{eq:theta-mu}.

For the dual statement, an $R$-linear map $g\colon L\to X$ is $\Lambda$-linear if and only if
\[
g(\theta(l'\otimes l))=f(l'\otimes g(l))
\qquad \forall l,l'\in L.
\]
This is precisely the kernel condition defining $\delta^\vee(A)=\Ker(\zeta_X-G(f^\sharp))$. Therefore $\delta^\vee(A)\cong\Hom_\Lambda(L,A)$ and the morphism $\gamma^\vee_A$ is exactly \eqref{eq:theta-nu}.
\end{proof}

For a class $\cC\subseteq R\text{-Mod}$, define
\begin{align*}
\Phi_\theta(\cC)
&=\{A\in\Lambda\text{-Mod}\mid \mu_A\text{ is monic and }A/LA\in\cC\},\\
\Psi_\theta(\cC)
&=\{A\in\Lambda\text{-Mod}\mid \mu^\vee_A\text{ is epic and }\operatorname{Ann}_A(L)\in\cC\}.
\end{align*}
By Proposition~\ref{prop:theta-functors}, we obtain that
\[
\Phi_\theta(\cC)=\Delta(\cC),
\qquad
\Psi_\theta(\cC)=\nabla(\cC),
\]
where the right-hand sides are taken in the extension and coextension descriptions of $\Lambda\text{-Mod}$ (see Section~\ref{section:extension} for the definitions of $\Delta(\cC)$ and $\nabla(\cC)$).
At the same time, we have the following description
\[
U^{-1}(\cC)=\{A\in \Lambda\Mod\mid A, \text{ viewed as an } R \text{-module via } R\hookrightarrow \Lambda, \text{ lies in }\cC\}.
\]

The four abstract cotorsion pairs of Section~\ref{section:extension} immediately acquire concrete module-theoretic descriptions.

\begin{corollary}\label{cor:theta-four-pairs}
Let $(\cX,\cY)$ be a cotorsion pair in $R\text{-Mod}$. Then 
\begin{enumerate}[label=\textup{(\arabic*)}]
\item $(\Phi_\theta(\cX),\Phi_\theta(\cX)^\perp)$ is a cotorsion pair in $\Lambda\text{-Mod}$.
\item $({}^\perp\Psi_\theta(\cY),\Psi_\theta(\cY))$ is a cotorsion pair in $\Lambda\text{-Mod}$.
\item If $\Tor_1^R(L,\cX)=0$, then
$({}^\perp U^{-1}(\cY),U^{-1}(\cY))$ is a cotorsion pair in $\Lambda\text{-Mod}$.
\item If $\Ext_R^1(L,\cY)=0$, then
$(U^{-1}(\cX),U^{-1}(\cX)^\perp)$  is a cotorsion pair in $\Lambda\text{-Mod}$.
\end{enumerate}
The cotorsion pairs in \textup{(3)} and \textup{(4)} are hereditary if and only if $(\cX,\cY)$ is hereditary. Moreover, if $(\cX,\cY)$ is hereditary, then
the cotorsion pair in \textup{(1)} is hereditary if and only if  $\Tor_1^\Lambda(L,\Phi_\theta(\cX))=0$ and the cotorsion pair in \textup{(2)} is hereditary if and only if $\Ext_\Lambda^1(L,\Psi_\theta(\cY))=0$.
\end{corollary}

\begin{proof}
Parts~(1) and~(2) follow from Theorems~\ref{thm:Delta-cotorsion} and \ref{thm:nabla-cotorsion}. Since $F=L\otimes_R-$, we have
$
\mathbb{L}_1F(X)=\Tor_1^R(L,X).
$
Thus part~(3) and its hereditary assertion are Theorem~\ref{thm:U-cotorsion}. Similarly, for $G=\Hom_R(L,-)$,
$
\mathbb{R}^1G(Y)=\Ext_R^1(L,Y),
$
so part~(4) and its hereditary assertion follow from Theorem~\ref{thm:coU-cotorsion}.

By Proposition~\ref{prop:theta-functors}, $\delta\cong L\otimes_\Lambda-$, so
$
\mathbb{L}_1\delta=\Tor_1^\Lambda(L,-).
$
The hereditary assertion for part~(1) is therefore Theorem~\ref{thm:Delta-heredity}. Dually, $\delta^\vee=\Hom_\Lambda(L,-)$, so
$
\mathbb{R}^1\delta^\vee=\Ext_\Lambda^1(L,-),
$
and Theorem~\ref{thm:nabla-heredity} gives the hereditary assertion for part~(2).
\end{proof}

\begin{remark}
{\rm (1) In the case $\theta=0$, Corollary~\ref{cor:theta-four-pairs} recovers \cite[Theorem~4.2]{MaoTrivial}.

(2) The Morita context ring $\Lambda_{(0,0)}=\begin{pmatrix}
A & {}_AN_B\\
{}_BM_A & B
\end{pmatrix}$ with zero bimodule homomorphisms is isomorphic to a $\theta$-extension ring $(A\times B)\ltimes_0 (M\oplus N)$; see
\cite[Section~7]{MaoTrivial}. In \cite{CuiRongZhang},  Zhang-Cui-Rong construct four cotorsion pairs over $\Lambda_{(0,0)}$. Via the above ring isomorphism, parts (3) and (4) of Corollary~\ref{cor:theta-four-pairs} recover
the construction of the two cotorsion pairs in \cite[Theorem~1.1]{CuiRongZhang}, while parts (1) and (2) of Corollary~\ref{cor:theta-four-pairs} generalize those in \cite[Theorem~1.2]{CuiRongZhang} by removing the condition
$M\otimes_AN=N\otimes_BM=0$.
}
\end{remark}

Let $\radf$ and $\coradf$ be the Beligiannis radical functor and coradical functor respectively. Set
\[
Q_n^\theta=C\radf^nT,
\qquad
(Q_n^\theta)^\vee=K\coradf^nH.
\]

\begin{proposition}\label{prop:theta-layers}
Let $A$ be a left $\Lambda$-module.
\begin{enumerate}[label=\textup{(\arabic*)}]
\item The Beligiannis radical functor is given by
$
\radf(A)\cong LA.
$
For $X\in R\text{-Mod}$,
\[
Q_n^\theta(X)
\cong
\frac{L^n(\Lambda\otimes_RX)}{L^{n+1}(\Lambda\otimes_RX)}.
\]

\item The Beligiannis coradical functor is given by
$
\coradf(A)\cong A/\Ann_A(L).
$
For $Y\in R\text{-Mod}$,
\[
(Q_n^\theta)^\vee(Y)
\cong
\frac{\Ann_{H(Y)}(L^{n+1})}{\Ann_{H(Y)}(L^n)}
\cong
\frac{\Hom_R(\Lambda/L^{n+1},Y)}
     {\Hom_R(\Lambda/L^n,Y)}.
\]
Here $\Hom_R(\Lambda/L^n,Y)$ is viewed as a submodule of $\Hom_R(\Lambda/L^{n+1},Y)$ via the injective map induced by the canonical epimorphism
\(\Lambda/L^{n+1}\twoheadrightarrow\Lambda/L^n\).
\end{enumerate}
\end{proposition}

\begin{proof}
Under the $\eta$-extension realization $(R\text{-Mod})\ltimes_\eta F$ of $\Lambda\text{-Mod}$, where
\(F=L\otimes_R-\), the structure map of a $\Lambda$-module $A$ is the multiplication map
$
L\otimes_RA\to A,
\, l\otimes a\mapsto la.
$
By definition, $\radf(A)$ is the image of this map; hence
\(\radf(A)=LA\) and $LA$ carries the induced $\Lambda$-module structure from $A$.
Therefore, for $X\in R\Mod$,
\[
Q_n^\theta(X)=C\radf^n(T(X))
=C(L^nT(X))
\cong
\frac{L^nT(X)}{L^{n+1}T(X)}=\frac{L^n(\Lambda\otimes_RX)}{L^{n+1}(\Lambda\otimes_RX)}.
\]

For the coextension realization, the structure map of $A$ is
\[
A\longrightarrow \Hom_R(L,A),
\qquad
a\longmapsto(l\mapsto la).
\]
The kernel is $K(A)=\Ann_A(L)$ as an $R$-module. Noting that $\Ann_A(L)$ can be viewed as a $\Lambda$-submodule of $A$  on which $L$ acts trivially since
$
l(\lambda a) =(l\lambda)a=0
$
for any $l\in L, \lambda\in \Lambda$ and $a\in \Ann_A(L)$, we have $ZK(A)=\Ann_A(L)$ as $\Lambda$-modules. Hence
\[
\coradf(A)=\Coker\bigl(ZK(A)\to A\bigr)
\cong A/\Ann_A(L).
\]
We shall prove $\coradf^n(A)\cong A/\Ann_A(L^n)$ by induction on $n$. The case $n=1$ is the preceding formula. Suppose that
$
\coradf^m(A)\cong A/\Ann_A(L^m).
$
An element $a+\Ann_A(L^m)$ belongs to
\(K\coradf^m(A)\) if and only if
$
L a\subseteq \Ann_A(L^m),
$
which is equivalent to
$
L^{m+1}a=0.
$
Therefore
\[
K\coradf^m(A)
\cong
\frac{\Ann_A(L^{m+1})}{\Ann_A(L^m)},
\]
Taking the corresponding quotient gives
\[
\begin{aligned}
\coradf^{m+1}(A)&=\Coker(ZK\coradf^m(A)\to \coradf^m(A))\\
&\cong
\frac{A/\Ann_A(L^m)}
     {\Ann_A(L^{m+1})/\Ann_A(L^m)}
\cong
A/\Ann_A(L^{m+1}),
\end{aligned}
\]
so the induction is complete. Then for $Y\in R\text{-Mod}$,
\[
(Q_n^\theta)^\vee(Y) = K\coradf^n(H(Y))
\cong
\frac{\Ann_{H(Y)}(L^{n+1})}{\Ann_{H(Y)}(L^n)}.
\]

For
\(\varphi\in H(Y)=\Hom_R(\Lambda,Y)\),
\[
\varphi\in\Ann_{H(Y)}(L^n)
\iff
l\varphi=0\text{ for every }l\in L^n
\iff
\varphi(\Lambda L^n)=0.
\]
Since $L^n$ is a two-sided ideal of the unital ring $\Lambda$, we have
\(\Lambda L^n=L^n\). Hence
\[
\Ann_{H(Y)}(L^n)
=
\{\varphi\in\Hom_R(\Lambda,Y)\mid \varphi(L^n)=0\}
\cong
\Hom_R(\Lambda/L^n,Y).
\]
Thus
\[
(Q_n^\theta)^\vee(Y)\cong \frac{\Hom_R(\Lambda/L^{n+1},Y)}{\Hom_R(\Lambda/L^n,Y)}.
\]
\end{proof}

We say that $\theta$ is \emph{$N$-nilpotent} if
\[
\theta\circ (\theta\otimes 1_L)\circ\cdots \circ (\theta\otimes 1_L\otimes \cdots\otimes 1_L)=0\colon L^{\otimes_R(N+1)}\to L.
\]
Equivalently, the corresponding natural transformation $\eta=\theta\otimes 1$ is $N$-nilpotent.

\begin{corollary}\label{cor:theta-completeness}
Assume that $\theta$ is $N$-nilpotent and that $(\cX,\cY)$ is a hereditary complete cotorsion pair in $R\text{-Mod}$.
\begin{enumerate}[label=\textup{(\arabic*)}]
\item If
$
\Tor_1^\Lambda(L,\Phi_\theta(\cX))=0,$ and
$
{L^n(\Lambda\otimes_RX)}/{L^{n+1}(\Lambda\otimes_RX)}\in\cY$ for all $X\in \cX$ and $1\le n\le N$,
then $(\Phi_\theta(\cX),\Phi_\theta(\cX)^\perp)$ is a hereditary complete cotorsion pair.
\item If
$
\Ext_\Lambda^1(L,\Psi_\theta(\cY))=0,$ and
${\Hom_R(\Lambda/L^{n+1},Y)}/{\Hom_R(\Lambda/L^n,Y)}$ $\in \cX$ for all $Y\in \cY$ and $1\le n\le N,$
then $({}^\perp\Psi_\theta(\cY),\Psi_\theta(\cY))$ is a hereditary complete cotorsion pair.
\end{enumerate}
\end{corollary}
\begin{proof}
Part~(1) is Theorem~\ref{thm:completeness} after using $\Phi_\theta(\cX)=\Delta(\cX)$, $\mathbb{L}_1\delta=\Tor_1^\Lambda(L,-)$ and $C\radf^nT\cong \frac{L^n(\Lambda\otimes_RX)}{L^{n+1}(\Lambda\otimes_RX)}$. Part~(2) is the corresponding dual statement from Theorem~\ref{thm:cocompleteness}, using $\Psi_\theta(\cX)=\nabla(\cX)$, $\mathbb{R}^1\delta^\vee=\Ext_\Lambda^1(L,-)$ and $K\coradf^nH\cong \frac{\Hom_R(\Lambda/L^{n+1},Y)}
     {\Hom_R(\Lambda/L^n,Y)}$.
\end{proof}

\begin{corollary}\label{cor:theta-U-completeness}
Assume that $(\cX,\cY)$ is a complete cotorsion pair in $R\text{-Mod}$.
\begin{enumerate}[label=\textup{(\arabic*)}]
\item If $\Tor_1^R(L,\cX)=0$ and
$
L\otimes_R\cX\subseteq\cY,
$
then $({}^\perp U^{-1}(\cY),U^{-1}(\cY))$ is a complete cotortion pair.
\item If $\Ext_R^1(L,\cY)=0$ and
$
\Hom_R(L,\cY)\subseteq\cX,
$
then $(U^{-1}(\cX),U^{-1}(\cX)^\perp)$ is a complete cotortion pair.
\end{enumerate}
\end{corollary}
\begin{proof}
Part~(1) is Theorem~\ref{thm:U-completeness}, because $\mathbb{L}_1F=\Tor_1^R(L,-)$ and $F(\cX)=L\otimes_R\cX$. Part~(2) is the corresponding dual statement from Theorem~\ref{thm:U-cocompleteness}, using $\mathbb{R}^1G=\Ext_R^1(L,-)$ and $G(\cY)=\Hom_R(L,\cY)$.
\end{proof}

By Theorem~\ref{thm:two-equal} and Theorem~\ref{thm:cotwo-equal}, we have the following results.

\begin{corollary}\label{cor:theta-pair-equality}
Let $(\cX,\cY)$ be a cotorsion pair in $R\text{-Mod}$. Assume that $\theta$ is nilpotent,
$
\Tor_1^R(L,\cX)=0,
$
and
$\Ext_R^1(A/LA,L\otimes_\Lambda A)=0$ for any $A\in\Phi_\theta(\cX).$
Then
$
(\Phi_\theta(\cX),\Phi_\theta(\cX)^\perp)
=({}^\perp U^{-1}(\cY),U^{-1}(\cY))
=(T(\cX),U^{-1}(\cY)).
$
\end{corollary}

\begin{corollary}
Let $(\cX,\cY)$ be a cotorsion pair in $R$-Mod. Assume that $\theta$ is nilpotent, $\Ext^1_R(L, \cY)=0$ and
$\Ext^1_R(\Hom_\Lambda(L, A), \Ann_A(L))=0$ for any $A\in \Psi_\theta(\cY)$. Then
$({}^\perp\Psi_\theta(\cY),\Psi_\theta(\cY))=(U^{-1}(\cX),U^{-1}(\cX)^\perp)=(U^{-1}(\cX),H(\cY))$.
\end{corollary}

\section{Cotorsion pairs over tensor rings}\label{section:tensor}

Let $R$ be a ring and $M$ an $R$-bimodule. The tensor ring of $M$ is defined by $T_R(M)=\bigoplus_{i\ge0}M^i$, where $M^{0}=R$ and $M^{i+1}=M\otimes_R M^{i}$ for $i\geq 0$.
The addition of elements in $T_R(M)$ is componentwise, and the multiplication is induced by the isomorphism
$(\bigoplus_{n=0}^\infty M^n)\otimes_R(\bigoplus_{m=0}^\infty M^m)\cong \bigoplus_{n,m=0}^\infty M^{n+m}$.
Set
$
M_+=\bigoplus_{i\ge1}M^i.
$
The multiplication gives an associative map
$
\theta\colon M_+\otimes_RM_+\longrightarrow M_+.
$
The tensor ring $\Lambda = T_R(M)$ can be viewed as a $\theta$-extension ring
$
R\ltimes_\theta M_+.
$
Thus Section~\ref{section:theta} applies with $L=M_+$. The category $\Lambda$-Mod is isomorphic to the $\eta$-extension $(R\Mod)\ltimes_\eta F$ with $F=M_+\otimes_R-
\cong\bigoplus_{i\ge1}M^i\otimes_R-,$ and $\eta=\theta\otimes 1$. $F$ has a right adjoint functor $G=\Hom_R(M_+,-)
\cong\prod_{i\ge1}\Hom_R(M^i,-).$

From \cite{ChenLu}, a left $\Lambda$-module can be written as a pair $\langle X,g\rangle$, where $X\in R\text{-Mod}$ and
$g\colon M\otimes_RX\longrightarrow X$
is $R$-linear. A morphism
$
\alpha\colon\langle X,g\rangle\longrightarrow\langle Y,h\rangle
$
is an $R$-linear map $\alpha\colon X\to Y$ satisfying
$
h\circ (1_M\otimes\alpha)=\alpha\circ g.
$

Given $g\colon M\otimes_RX\to X$, define $g_i\colon M^i\otimes_RX\to X$ recursively by $g_0=1_X$ and
$g_i=g\circ(1_M\otimes g_{i-1}).$
Here we identify $M^i\otimes_R X$ with $M\otimes_R M^{i-1}\otimes_R X$. Set
\[
g^e=(g_1,g_2,\cdots)
\colon M_+\otimes_RX=\bigoplus_{i\geq 1}M^i\otimes_RX\longrightarrow X.
\]
Then the object $\langle X,g\rangle \in \Lambda\Mod$ is identified with $(X,g^e)\in (R\Mod)\ltimes_\eta F$.
Additionally, define
\[
g^\sharp\colon X\longrightarrow\Hom_R(M,X),
\qquad
g^\sharp(x)(m)=g(m\otimes x).
\]

\begin{proposition}\label{prop:tensor-bimodule}
Let $\Lambda=T_R(M)$ be a tensor ring. For every $A=\langle X,g\rangle\in\Lambda\text{-Mod}$,
\begin{align*}
C(A)=A/M_+A\cong\Coker g,\qquad &\delta(A)\cong M_+\otimes_\Lambda A\cong M\otimes_RX,\\
K(A)=\Ann_A(M_+)\cong\Ker g^\sharp,\qquad &\delta^\vee(A)\cong\Hom_\Lambda(M_+,A)\cong\Hom_R(M,X).
\end{align*}
Under these isomorphisms the multiplication map $\mu_A$ and its adjoint $\mu_A^\vee$ in Proposition~\ref{prop:theta-functors} become $g$ and $g^\sharp$, respectively.
\end{proposition}

\begin{proof}
Since every $g_i$ factors through $g$, we have $M_+A=\Ima g$. Moreover, an element is annihilated by $M_+$ if and only if it is annihilated by $M$. Hence
\[
C(A)=A/M_+A\cong\Coker g,
\qquad
K(A)=\Ann_A(M_+)\cong\Ker g^\sharp.
\]

There are natural bimodule isomorphisms
\[
M\otimes_R\Lambda\xrightarrow{\ \sim\ }M_+,
\qquad m\otimes\lambda\longmapsto m\lambda,
\]
as $R$--$\Lambda$-bimodules, and
\[
\Lambda\otimes_RM\xrightarrow{\ \sim\ }M_+,
\qquad \lambda\otimes m\longmapsto\lambda m,
\]
as $\Lambda$--$R$-bimodules.
Tensoring the first isomorphism over $\Lambda$ with $A$ gives
\[
M_+\otimes_\Lambda A
\cong M\otimes_R\Lambda\otimes_\Lambda A
\cong M\otimes_RX,
\]
and applying $\Hom_\Lambda(-,A)$ to the second one gives
\[
\Hom_\Lambda(M_+,A)
\cong\Hom_\Lambda(\Lambda\otimes_RM,A)
\cong\Hom_R(M,X).
\]
The multiplication map $\mu_A$ sends $m\otimes x$ to $g(m\otimes x)$, so it becomes $g$, and its adjoint $\mu_A^\vee$ becomes $g^\sharp$.
\end{proof}

For a class $\cC\subseteq R\text{-Mod}$, define
\begin{align*}
\Phi(\cC)
&=\{\langle X,g\rangle\mid g\text{ is monic and }\Coker g\in\cC\},\\
\Psi(\cC)
&=\{\langle X,g\rangle\mid g^\sharp\text{ is epic and }\Ker g^\sharp\in\cC\}.
\end{align*}
By Proposition~\ref{prop:tensor-bimodule}, we have $\Phi(\cC)=\Phi_\theta(\cC)$ and $\Psi(\cC)=\Psi_\theta(\cC)$. This naturally yields the following cotorsion pairs over tensor rings.

\begin{corollary}\label{thm:tensor-four-pairs}
Let $\Lambda=T_R(M)$ be a tensor ring and $(\cX,\cY)$ a cotorsion pair in $R\text{-Mod}$. Then 
\begin{enumerate}[label=\textup{(\arabic*)}]
\item $(\Phi(\cX),\Phi(\cX)^\perp)$ is a cotorsion pair in $\Lambda\text{-Mod}$.
\item $({}^\perp\Psi(\cY),\Psi(\cY))$ is a cotorsion pair in $\Lambda\text{-Mod}$.
\item If
$
\Tor_1^R(M^i,\cX)=0$ for all $i\ge 1$,
then $({}^\perp U^{-1}(\cY),U^{-1}(\cY))$ is a cotorsion pair in $\Lambda\text{-Mod}$.
\item If
$
\Ext_R^1(M^i,\cY)=0$ for all $i\ge 1$,
then $(U^{-1}(\cX),U^{-1}(\cX)^\perp)$ is a cotorsion pair in $\Lambda\text{-Mod}$.
\end{enumerate}
The cotorsion pairs in \textup{(3)} and \textup{(4)} are hereditary if and only if $(\cX,\cY)$ is hereditary. The condition in \textup{(3)} holds if $M_R$ is flat and the condition in \textup{(4)} holds if ${}_RM$ is projective.
\end{corollary}

\begin{proof}
Note that $\Tor_1^R(M_+,X)
\cong\bigoplus_{i\ge1}\Tor_1^R(M^i,X),
$ and $
\Ext_R^1(M_+,Y)
\cong\prod_{i\ge1}\Ext_R^1(M^i,Y)
$. Then these results follow from Corollary~\ref{cor:theta-four-pairs}.
\end{proof}

\begin{corollary}\label{cor:tensor-heredity}
Let $\Lambda=T_R(M)$ be a tensor ring. Assume that $(\cX,\cY)$ is a hereditary cotorsion pair.
\begin{enumerate}[label=\textup{(\arabic*)}]
\item The cotorsion pair $(\Phi(\cX),\Phi(\cX)^\perp)$ is hereditary if and only if
$
\Tor_1^\Lambda(M_+,\Phi(\cX))=0.
$
In particular, it is hereditary if $\Tor_1^R(M,U(\Phi(\cX)))=0$.
\item The cotorsion pair $({}^\perp\Psi(\cY),\Psi(\cY))$ is hereditary if and only if
$
\Ext_\Lambda^1(M_+,\Psi(\cY))=0.
$
In particular, it is hereditary if $\Ext^1_R(M,U(\Psi(\cY)))=0$.
\end{enumerate}
\end{corollary}

\begin{proof}
The general statements follow from the hereditary assertions in Corollary~\ref{cor:theta-four-pairs} with $L=M_+$. Assume that
$\Tor_1^R(M,U(\Phi(\cX)))=0$. Let $0\to A_1\to A_2\to A\to 0$ be a short exact sequence of $\Lambda$-modules with $A\in \Phi(\cX)$.
Applying the exact functor $U$, we get a short exact sequence $0\to U(A_1)\to U(A_2)\to U(A)\to 0$. Since $\Tor_1^R(M,U(A))=0$, the
sequence $0\to M\otimes_RU(A_1)\to M\otimes_RU(A_2)\to M\otimes_RU(A)\to 0$ is still exact.
Note that
$
M_+\otimes_\Lambda-
\cong M\otimes_RU(-).
$
It follows that $\Tor_1^\Lambda(M_+,\Phi(\cX))=0$ from Lemma~\ref{lem:L_1-seq} and hence $(\Phi(\cX),\Phi(\cX)^\perp)$ is hereditary. Dually, since
$
\Hom_\Lambda(M_+,-)
\cong\Hom_R(M,U(-))
$, we obtain the sufficient condition in (2).
\end{proof}

\begin{proposition}\label{prop:tensor-layers}
Let $\Lambda=T_R(M)$ be a tensor ring. For $X,Y\in R\text{-Mod}$ and $n\geq 1$, there are natural isomorphisms
\begin{equation*}
\begin{aligned}
Q_n^\theta(X)&=C\radf^nT(X)\cong M^n\otimes_RX,\\
(Q_n^\theta)^\vee(Y)&=K\coradf^nH(Y)\cong\Hom_R(M^n,Y).
\end{aligned}
\end{equation*}
\end{proposition}

\begin{proof}
Applying Proposition~\ref{prop:theta-layers} with $L=M_+$, we have
\[
Q_n^\theta(X)\cong\frac{M_+^n(\Lambda\otimes_RX)}{M_+^{n+1}(\Lambda\otimes_RX)},\quad (Q_n^\theta)^\vee(Y)
\cong
\frac{\Hom_R(\Lambda/M_+^{n+1},Y)}
     {\Hom_R(\Lambda/M_+^n,Y)}.
\]
Since
$
M_+^n=\bigoplus_{i=n}^\infty M^i,
$
and
$
M_+^{n+1}=\bigoplus_{i=n+1}^\infty M^i,
$
both $M_+^n$ and $M_+^{n+1}$ are direct summands of $\Lambda$ as right $R$-modules. Hence their inclusions into $\Lambda$ remain monic after tensoring with $X$. It follows that
\begin{align*}
&M_+^n(\Lambda\otimes_RX)\\
=&\Ima(M_+^n\otimes_R(\Lambda\otimes_RX)\to \Lambda\otimes_RX)\\
=&\Ima(M_+^n\otimes_RX\to \Lambda\otimes_RX)\cong M_+^n\otimes_RX.
\end{align*}
Similarly,
$
M_+^{n+1}(\Lambda\otimes_RX)\cong M_+^{n+1}\otimes_RX.
$
Moreover, the split exact sequence of right $R$-modules
\[
0\longrightarrow M_+^{n+1}\longrightarrow M_+^n
\longrightarrow M_+^n/M_+^{n+1}\longrightarrow0
\]
remains exact after tensoring with $X$. Therefore
\[
\frac{M_+^n(\Lambda\otimes_RX)}{M_+^{n+1}(\Lambda\otimes_RX)}
\cong
\frac{M_+^n\otimes_RX}{M_+^{n+1}\otimes_RX}
\cong
\left(M_+^n/M_+^{n+1}\right)\otimes_RX
\cong M^n\otimes_RX,
\]
which gives the first formula.

We now prove the dual formula. Since
$
\Lambda/M_+^{n+1}\cong\bigoplus_{i=0}^nM^i,
$
and
$
\Lambda/M_+^n\cong\bigoplus_{i=0}^{n-1}M^i,
$
we have
\[
(Q_n^\theta)^\vee(Y)
\cong
\frac{\Hom_R(\Lambda/M_+^{n+1},Y)}
     {\Hom_R(\Lambda/M_+^n,Y)}
\cong\Hom_R(M^n,Y).
\]
\end{proof}

If
$
M^{N+1}=0
$
for some $N\ge1$, we say that $M$ is \emph{$N$-nilpotent}. Then
$
M_+=\bigoplus_{i=1}^NM^i,
$
and the map $\theta\colon M_+\otimes_RM_+\to M_+$ is $N$-nilpotent.

From Corollary~\ref{cor:theta-completeness}, Proposition~\ref{prop:tensor-layers} and Corollary~\ref{cor:theta-U-completeness}, the following results are immediate.

\begin{corollary}\label{thm:tensor-completeness}
Let $\Lambda = T_R(M)$ be a tensor ring with $M$ being $N$-nilpotent. Assume that $(\cX,\cY)$ is a hereditary complete cotorsion pair in $R\text{-Mod}$.
\begin{enumerate}[label=\textup{(\arabic*)}]
\item If
$
\Tor_1^\Lambda(M_+,\Phi(\cX))=0
$
and
$
M^i\otimes_R\cX\subseteq\cY$ for all $1\le i\le N,$
then $(\Phi(\cX),\Phi(\cX)^\perp)$ is a hereditary complete cotorsion pair.
\item If
$
\Ext_\Lambda^1(M_+,\Psi(\cY))=0
$
and
$
\Hom_R(M^i,\cY)\subseteq\cX$ for all $1\le i\le N,$
then $({}^\perp\Psi(\cY),\Psi(\cY))$ is a hereditary complete cotorsion pair.
\end{enumerate}
\end{corollary}

\begin{corollary}\label{thm:tensor-U-completeness}
Let $\Lambda = T_R(M)$ be a tensor ring. Assume that $(\cX,\cY)$ is a complete cotorsion pair in $R\text{-Mod}$.
\begin{enumerate}[label=\textup{(\arabic*)}]
\item If
$
\Tor_1^R(M^i,\cX)=0
$
and
$
M^i\otimes_R\cX\subseteq\cY
$ for all $i\ge 1$,
then $({}^\perp U^{-1}(\cY),U^{-1}(\cY))$ is a complete cotorsion pair.
\item If
$
\Ext_R^1(M^i,\cY)=0
$
and
$
\Hom_R(M^i,\cY)\subseteq\cX
$ for all $i\ge 1$,
then $(U^{-1}(\cX),U^{-1}(\cX)^\perp)$ is a complete cotorsion pair.
\end{enumerate}
\end{corollary}

Now we consider the conditions under which the two cotorsion pairs $(\Phi(\cX),\Phi(\cX)^\perp)$ and $({}^\perp U^{-1}(\cY), U^{-1}(\cY))$ coincide.
Let $A=\langle E,g\rangle\in\Phi(\cC)$ and $X=\Coker g$. By Propositions~\ref{prop:tensor-bimodule} and \ref{prop:theta-functors}, the sequence \eqref{eq:T-canonical seq} becomes
\begin{equation}\label{eq:tensor-basic}
0\longrightarrow M\otimes_RE
\xlongrightarrow{g}E
\longrightarrow X\longrightarrow0,
\end{equation}
Consequently, by Proposition~\ref{prop:Delta-induced}, $\Phi(\cC)=T(\cC)$ exactly when this canonical sequence \eqref{eq:tensor-basic}
is split exact for every $A\in\Phi(\cC)$. The next proposition gives convenient homological conditions forcing the required splitting.

\begin{proposition}\label{prop:tensor-Phi-T}
Let $\Lambda = T_R(M)$ be a tensor ring with $M$ being $N$-nilpotent, and let $\cC$ be a class in $R\text{-Mod}$. Assume that for every $X\in\cC$,
$\Tor_1^R(M^i,X)=0$ for all $1\le i\le N-1$,
and
$
\Ext_R^1(X,M^i\otimes_RX)=0$ for all $1\le i\le N$.
Then
$
\Phi(\cC)=T(\cC).
$
\end{proposition}

\begin{proof}
Let $A=\langle E,g\rangle\in\Phi(\cC)$ and $X=\Coker g\in\cC$. Since $g$ is monic, we have the exact sequence \eqref{eq:tensor-basic}.

For $1\le i\le N$, set
$
K_i=M^i\otimes_RE.
$
For $1\le i\le N-1$, applying $M^i\otimes_R-$ to \eqref{eq:tensor-basic} and using $\Tor_1^R(M^i,X)=0$ gives a short exact sequence
\begin{equation}\label{eq:tensor-Ki}
0\longrightarrow K_{i+1}\longrightarrow K_i
\longrightarrow M^i\otimes_RX\longrightarrow0.
\end{equation}
Since $M^{N+1}=0$, applying $M^N\otimes_R-$ to \eqref{eq:tensor-basic} shows that
$
K_N\cong M^N\otimes_RX.
$
Hence $\Ext_R^1(X,K_N)=0$. We shall prove that $\Ext_R^1(X,K_i)=0$ for $1\leq i\leq N$ by descending induction on $i$.
Suppose $\Ext_R^1(X,K_{i+1})=0$. Applying $\Hom_R(X,-)$ to \eqref{eq:tensor-Ki} yields an exact sequence
\[
\Ext_R^1(X,K_{i+1})\longrightarrow
\Ext_R^1(X,K_i)\longrightarrow
\Ext_R^1(X,M^i\otimes_RX).
\]
The two outer terms vanish, so $\Ext_R^1(X,K_i)=0$. Therefore
\[
\Ext_R^1(X,M\otimes_RE)=\Ext_R^1(X,K_1)=0.
\]
The sequence \eqref{eq:tensor-basic} splits. Hence $\Phi(\cC)=T(\cC).$
\end{proof}

\begin{corollary}\label{cor:tensor-coincidence}
Let $\Lambda = T_R(M)$ be a tensor ring with $M$ being $N$-nilpotent, and let $(\cX,\cY)$ be a cotorsion pair in $R\text{-Mod}$. Assume that
$
\Tor_1^R(M^i,X)=0
$
and
$
\Ext_R^1(X,M^i\otimes_RX)=0
$ for all $X\in\cX$ and $1\le i\le N$.
Then
\[
(\Phi(\cX),\Phi(\cX)^\perp)
=({}^\perp U^{-1}(\cY),U^{-1}(\cY))
=(T(\cX),U^{-1}(\cY)).
\]
\end{corollary}

\begin{proof}
By Proposition~\ref{prop:tensor-Phi-T} we have $\Phi(\cX)=T(\cX)$. The Tor-vanishing assumptions give
$
\mathbb{L}_1F(\cX)=\Tor_1^R(M_+,\cX)=0.
$
Hence the proof of Proposition~\ref{prop:relation} yields the equality of the two cotorsion pairs.
\end{proof}

Dually, we have the following.

\begin{corollary}\label{thm:tensor-cocoincidence}
Let $\Lambda = T_R(M)$ be a tensor ring with $M$ being $N$-nilpotent, and let $(\cX,\cY)$ be a cotorsion pair in $R\text{-Mod}$. Assume that
$
\Ext^1_R(M^i,Y)=0
$
and
$
\Ext_R^1(\Hom_R(M^i,Y),Y)=0
$ for all $Y\in\cY$ and $1\le i\le N$.
Then
\[
({}^\perp\Psi(\cY),\Psi(\cY))
=(U^{-1}(\cX),U^{-1}(\cX)^\perp)
=(U^{-1}(\cX),H(\cY)).
\]
\end{corollary}

\section*{Acknowledgments}
This work is supported by
the National Natural Science Foundation of China (12561008),
the project of Young and Middle-aged Academic and Technological leader of Yunnan
(202305AC160005) and the Basic Research Program of Yunnan Province (202301AT070070).

\end{document}